\documentclass{elsart}

\usepackage{amsfonts}
\usepackage{amsmath}
\usepackage{amssymb}
\usepackage[all,dvips]{xy}

\newcommand{\DCC}{{\ensuremath{\Ext^2_C(DC,C)}}}
\newcommand{\Ctilde}{\ensuremath{\tilde{C}}}
\newcommand{\Cbar}{{\ensuremath{\overline{C}}}}
\newcommand{\Pbar}{{\ensuremath{\overline{P}}}}

\newcommand{\Crep}{{\ensuremath{\stackrel{\scriptscriptstyle\vee}{C}}}}
\newcommand{\Prep}{{\ensuremath{\stackrel{\scriptscriptstyle\vee}{P}}}}
\newcommand{\Irep}{{\ensuremath{\stackrel{\scriptscriptstyle\vee}{I}}}}
\newcommand{\dc}{{\ensuremath{\stackrel{\scriptscriptstyle\vee}{d}}}}

\newcommand{\rad}{\ensuremath{\textup{rad}\,}}
\newcommand{\Hom}{\ensuremath{\textup{Hom}\,}}
\newcommand{\sHom}{\ensuremath{\underline{\textup{Hom}}\,}}
\newcommand{\Ext}{\ensuremath{\textup{Ext}}}
\newcommand{\End}{\ensuremath{\textup{End}\,}}
\newcommand{\Ann}{\ensuremath{\textup{Ann}\,}}
\newcommand{\Supp}{\ensuremath{\textup{Supp}\,}}

\newcommand{\ind}{\ensuremath{\textup{ind}\,}}
\newcommand{\add}{\ensuremath{\textup{add}\,}}
\newcommand{\iadd}{\ensuremath{\textup{iadd}\,}}

\newcommand{\pd}{\ensuremath{\textup{pd}\,}}

\newcommand{\CA}{{\ensuremath{\mathcal{C}_A}}}
\newcommand{\CC}{{\ensuremath{\mathcal{C}_C}}}
\newcommand{\DA}{{\ensuremath{\mathcal{D}^b(\textup{mod}\, A)}}}
\newcommand{\DC}{{\ensuremath{\mathcal{D}^b(\textup{mod}\, C)}}}
\newcommand{\DCprime}{{\ensuremath{\mathcal{D}^b(\textup{mod}\, C')}}}
\newcommand{\DkS}{{\ensuremath{\mathcal{D}^b(\textup{mod}\, k\Sigma)}}}

\newcommand{\dcc}{\Ext^2_C(DC,C)}
\newcommand{\cdcc}{C\ltimes\Ext^2_C(DC,C)}

\newcommand{\zG}{\ensuremath{\Gamma}}
\newcommand{\zS}{\ensuremath{\Sigma}}
\newcommand{\zO}{\ensuremath{\Omega}}
\newcommand{\za}{\ensuremath{\alpha}}
\newcommand{\zb}{\ensuremath{\beta}}
\newcommand{\ze}{\ensuremath{\epsilon}}
\newcommand{\zd}{\ensuremath{\delta}}
\newcommand{\zg}{\ensuremath{\gamma}}
\newcommand{\zl}{\ensuremath{\lambda}}

\newcommand{\zs}{\ensuremath{\sigma}}

\newcommand{\Sp}{\ensuremath{\Sigma^+}}
\newcommand{\sx}{\ensuremath{\sigma^+_x}}

\newcommand{\sC}{\ensuremath{\sigma^+_xC}}

\begin{document}

\begin{frontmatter}

\title{Cluster-tilted algebras without clusters}
  \author{Ibrahim Assem\thanksref{nserc}}, \address{D\'epartement de Math\'ematiques,
  Universit\'e de Sherbrooke, Sherbrooke (Qu\'ebec), J1K 2R1, Canada}
  \ead{ibrahim.assem@usherbrooke.ca}
  \author{Thomas Br\"ustle\thanksref{nserc+}}, \address{D\'epartement de Math\'ematiques,
  Universit\'e de Sherbrooke, Sherbrooke (Qu\'ebec), J1K 2R1, Canada \emph{and} 
  Department of Mathematics, Bishop's University,  Lennoxville, (Qu\'ebec),
  J1M 1Z7, Canada} 
  \ead{thomas.brustle@usherbrooke.ca}
  \author{Ralf Schiffler\thanksref{uconn}}, \address{Department of Mathematics, University of Connecticut, Storrs CT, 06269-3009, USA}
  \ead{schiffler@math.uconn.edu} 

  \thanks[nserc]{Partially supported by the NSERC of Canada and the Universit\'e de Sherbrooke}
  \thanks[nserc+]{Partially supported by the NSERC of Canada and the universities of Sherbrooke and Bishop's}
  \thanks[uconn]{Partially supported by the NSF grant DMS-0908765 and the University of Connecticut}

\begin{abstract} 
\end{abstract}
\end{frontmatter}


\section{Introduction}\label{sect 0}

Cluster-tilted algebras were introduced in \cite{BMR} and also, independently, in \cite{CCS} for type $\mathbb{A}$, as a by-product of the theory of cluster algebras of Fomin and Zelevinsky \cite{FZ}. They are the endomorphism algebras of the so-called tilting objects in the cluster category of \cite{BMRRT}. Since their introduction, they have been the subject of several investigations, see, for instance, \cite{BMR,CCS, ABS, KR, BFPPT,BOW}. Part of their interest comes from the fact that the cluster category is a 2-Calabi-Yau category. In particular, the representation theory of cluster-tilted algebras has been shown to be very similar to that of the self-injective algebras, see \cite{ABS,ABS2, ABS3}. One of the essential tools in the study of self-injective algebras is the notion of reflection of a tilted algebra, introduced by Hughes and Waschb\"usch in \cite{HW}. This allowed to prove that, if $C$ is a tilted algebra, then its trivial extension $T(C)$ by the minimal injective cogenerator  bimodule is representation-finite if and only if $C$ is of Dynkin type and, in this case, $T(C)\cong T(B)$ if and only if $B$ is an iterated reflection of $C$ (or, equivalently, $B$ is iterated tilted of the same type as $C$), see also \cite{BLR, AHR, Ho}. Moreover, the proofs of these results developed into algorithms allowing to compute explicitly the module category of $T(C)$, starting from that of $C$, see \cite{HW, BLR}.

We recall from \cite{ABS} that, if $C$ is a tilted algebra, then the trivial extension $\Ctilde$ of $C$ by the $C$-$C$-bimodule $\dcc$ is cluster-tilted, and conversely, every cluster-tilted algebra is of this form. 
On the other hand, this (surjective) map from tilted algebras to cluster-tilted algebras is certainly not injective and it is an interesting question to find all the tilted algebras $B$ such that $\tilde B = \Ctilde.$ This  problem has already been considered in \cite{ABS2} and \cite{BOW}, see also \cite{BFPPT}. In the present paper, we define notions of  reflections (and, dually coreflections) of complete slices and of tilted algebras. Our main result may now be stated as follows.

\begin{thm}
Let  $C$ be a tilted algebra having a tree $\zS$ as a complete slice. A tilted algebra $B$ is such that $\tilde B=\Ctilde$ if and only if there exists a sequence of reflections and coreflections $\zs_1,\ldots,\zs_t$ such that $B=\zs_1\cdots\zs_t C$ has $\zO=\zs_1\cdots\zs_t\zS$ as a complete slice and $B=\Ctilde/\Ann \zO$.
\end{thm}

The restriction to tilted algebras of tree type seems to be necessary to ensure the existence of reflections.

As a consequence of this construction and our proof, we obtain, as in \cite{HW}, an algorithm allowing to compute explicitly the transjective component of the module category of $\Ctilde$, having as starting data only the knowledge of the tilted algebra $C$. In particular, if $C$ is of Dynkin type, this yields the whole module category of $\Ctilde$. We observe that, since the transjective component of the module category of $\Ctilde$ is standard, then it is uniquely determined by combinatorial data.

The paper is organised as follows. After a short preliminary section, in which we fix the notation and recall the needed results, we devote our section \ref{sect 2} to general properties of the Auslander-Reiten quiver of a cluster-tilted algebra. In section \ref{sect 3}, we define reflections of complete slices and of tilted algebras. Section \ref{sect 4} is devoted to the proof of our main results, and section \ref{sect 5} to the algorithm. We end the paper in section \ref{sect 6} by showing how our algorithm may be applied to construct the tubes of cluster-tilted algebras of euclidean type.


\section{Preliminaries}\label{sect 1} \begin{subsection}{Notation}\label{sect 1.1}
Throughout this paper, algebras are basic and connected, locally finite dimensional over an algebraically closed field $k$. For an algebra $C$, we denote by $\textup{mod}\,C$ the category of finitely generated right $C$-modules. All subcategories are full and so are identified with their object classes. Given a category $\mathcal{C}$, we sometimes write $M\in \mathcal{C}$ to express that $M$ is an object in $\mathcal{C}$. If $\mathcal{C}$ is a full subcategory of $\textup{mod}\,C$, we denote by $\add\mathcal{C}$ the full subcategory of $\textup{mod}\,C$ having as objects the finite direct sums of summands of modules in $\mathcal{C}$.

Following \cite{BoG}, we sometimes consider equivalently an algebra $C$ as a locally bounded $k$-category, in which the object class $C_0$ is (in bijection with) a complete set $\{e_x\}$ of primitive orthogonal idempotents of $C$, and the space of morphisms from $e_x$ to $e_y$ is $C(x,y)=e_xCe_y$. A full subcategory $B$ of $C$ is \emph{convex} if, for any path $x=x_0\to x_1\to \ldots \to x_t=y$ in the quiver $Q_C$ of $C$, with $x,y\in B$, we have $x_i\in B$ for all $i$. For a point $x$ in $Q_C$, we denote by $P_x,I_x,S_x$ respectively the indecomposable projective, injective and simple $C$-modules corresponding to $x$. We denote by $\zG(\textup{mod}\,C)$ the Auslander-Reiten quiver of $C$ and by $\tau_C=D\,Tr, \tau_C^{-1}=Tr\,D$ the Auslander-Reiten translations. Given two points $M,N$ in $\zG(\textup{mod}\,C)$, we denote by $M \leadsto N$ or by $M\le N$ the existence of a path (of non-zero morphisms between indecomposable modules) from $M$ to $N$ in $\textup{mod}\,C$. More generally, if $\mathcal{S}_1,\mathcal{S}_2$ are two sets of indecomposable modules, we write $\mathcal{S}_1\le \mathcal{S}_2$ if every module in $\mathcal{S}_1$ has a successor in $\mathcal{S}_2$, no module in $\mathcal{S}_2$ has a successor in $\mathcal{S}_1$, and no module in $\mathcal{S}_1$ has a predecessor in $\mathcal{S}_2$. The notation $\mathcal{S}_1<\mathcal{S}_2$ stands for $\mathcal{S}_1\le \mathcal{S}_2$ and $\mathcal{S}_1\cap \mathcal{S}_2=\emptyset$.

For further definitions and facts, we refer the reader to \cite{ARS,ASS}. For tilting theory, we refer to \cite{ASS, Ri}.

\end{subsection}

\begin{subsection}{Cluster-tilted algebras}\label{sect 1.2} Let $A$ be a finite dimensional hereditary $k$-algebra, The \emph{cluster category} $\CA$ of $A$ is defined as follows. Let $F$ be the automorphism of the bounded derived category $\DA$ defined as the composition $\tau^{-1}_{\mathcal{D}}[1]$, where $\tau^{-1}_{\mathcal{D}}$ is the Auslander-Reiten translation in $\DA$ and $[1]$ is the shift (suspension) functor. Then $\CA$ is the orbit category $\DA/F$, its objects are the $F$-orbits $\tilde X=(F^i X)_{i\in \mathbb{Z}} $ of the objects $X\in \DA$ and the space of morphisms from $\tilde X=(F^iX)_i$ to   $\tilde Y=(F^iY)_i$ is
\[\Hom_{\CA}(\tilde X,\tilde Y)= \oplus_{i\in\mathbb{Z}} \Hom_{\DA}(X,F^iY).\]
$\CA$ is a triangulated Krull-Schmidt category with almost split triangles. The projection $\pi:\DA\to \CA$ is a triangle functor which commutes with the Auslander-Reiten translations \cite{BMRRT,K}. Moreover, for any two objects $\tilde X,\tilde Y$ in $\CA$, we have a functorial isomorphism $\Ext^1_{\CA}(\tilde X,\tilde Y)\cong D\Ext^1_{\CA}(\tilde Y,\tilde X)$, in other words, the category $\CA$ is 2-Calabi-Yau.

An object $\tilde T\in \CA$ is \emph{tilting} if $\Ext^1_{\CA}(\tilde T,\tilde T)=0$, and the number of isomorphism classes of indecomposable summands of $\tilde T$ equals the rank of the Grothendieck group $K_0(A)$ of $A$. The endomorphism algebra $B=\End_{\CA} \tilde T$ is then called \emph{cluster-tilted}. Moreover, we have an equivalence $\textup{mod}\,B\cong \CA/\iadd (\tau_{\CA}\tilde T)$, where $\tau_{\CA}$ is the Auslander-Reiten translation in $\CA$ and $\iadd (\tau_{\CA}\tilde T)$ is the ideal of $\CA$ consisting of all morphisms factoring through objects of $\add (\tau_{\CA}\tilde T)$. Also, this equivalence commutes with the Auslander-Reiten translations in both categories \cite{BMR}.

We now describe the Auslander-Reiten quivers of $\CA$ and $B$. If $A=kQ$ is representation-finite, the $\zG(\CA)$ is of the form $\mathbb{Z}Q/\langle\varphi\rangle$, where $\varphi$ is the automorphism of $\mathbb{Z}Q$ induced by $F$. If $A=kQ$ is representation infinite, then $\zG(\CA)$ has a unique component of the form $\mathbb{Z}Q$, called \emph{transjective}, because it is the image (under $\pi$) of the transjective components of $\zG(\DA)$. Moreover,  $\zG(\CA)$ also has components called \emph{regular}, because they are the image of the regular components of  $\zG(\CA)$. In both cases, we deduce $\zG(\textup{mod}\,B)$ from  $\zG(\CA)$ by simply deleting the $|Q_0|$ points corresponding to the summands of  $\tau_{\CA}\tilde T$.  
\end{subsection}

\begin{subsection}
{Relation-extensions and slices}\label{sect 1.3}
If $B$ is cluster-tilted, then there exists a hereditary algebra $A$ and a tilting $A$-module $T$ such that $B=\End_{\CA} \tilde T$, see \cite[3.3]{BMRRT}. Moreover, if $C=\End_A T$ is the corresponding tilted algebra, then the trivial extension $\tilde C=\cdcc$ (the \emph{relation-extension} of $C$) is cluster-tilted and, actually, isomorphic to $B$, see \cite{ABS}. Now, tilted algebras are characterised by the  presence of so-called complete slices in the connecting components of their Auslander-Reiten quivers \cite{ASS, Ri}. The corresponding notion for cluster-tilted algebras is as follows \cite{ABS2}. A full subquiver $\zS$ of $\zG(\textup{mod}\,\Ctilde) $ is a \emph{local slice} if :
\begin{itemize}
\item[(LS1)] $\zS $ is a presection, that is \begin{itemize}
\item[(a)] If $X\in \zS$ and $X\to Y$ is an arrow, then either $Y\in \zS$ or $\tau_{\Ctilde}Y\in \zS$.
\item[(b)] If $Y\in \zS$ and $X\to Y$ is an arrow, then either $X\in \zS$ or $\tau^{-1}_{\Ctilde}X\in \zS$.
\end{itemize}
\item[(LS2)] $\zS$ is sectionally convex, that is, if $X=X_0\to X_1\to\cdots\to X_t=Y$ is a sectional path in $\zG(\textup{mod}\,\Ctilde)$, with $X,Y\in \zS$, then $X_i \in  \zS$ for all $i$.
\item[(LS3)] $|\zS_0|=\textup{rk} K_0(C)$.
\end{itemize}
Let $C$ be tilted, then, under the standard embedding $\textup{mod}\,C\to \textup{mod}\,\Ctilde$ any complete slice in $\textup{mod}\,C$ embeds as a local slice in $\textup{mod}\,\Ctilde$, and any local slice occurs in this way. If $B$ is cluster-tilted, then a tilted algebra $C$ is such that $B=\Ctilde$ if and only if there exists a local slice $\zS$ in $\zG(\textup{mod}\,B)$ such that $C=B/\Ann\!_B\zS$, where $\Ann\!_B\zS =\cap_{X\in\zS} \Ann\!_B X$, see \cite{ABS2}.

\end{subsection}

\begin{subsection}{Cluster-repetitive algebras}\label{sect 1.3} Let $C$ be a tilted algebra. Its  \emph{cluster-repetitive algebra} $\Crep$ 
 is the locally finite dimensional algebra given by
\[
\Crep  \quad = \quad \left[
  \begin{array}{cccccccccc}
  \ddots &&&\ 0\ \\
&\ C_{-1}\ \\
&E_0&\ C_0\ \\
&&E_{1}&\ C_1\ \\
&\ 0\ &&&\ddots
  \end{array}
 \right] \]
where matrices have only finitely many non-zero coefficients, $C_i=C$ 
and $E_i =\Ext^2_C(DC,C) $ for all $i\in \mathbb{Z}$, all the remaining
coefficients are zero, and the multiplication is induced from that of $C$,
the $C$-$C$-bimodule structure of $\Ext^2_C(DC,C) $ and the zero map
$\Ext^2_C(DC,C) \otimes_C \Ext^2_C(DC,C) \to 0$.
The identity maps $C_i\to C_{i-1}$, $E_i\to E_{i-1}$ induce an automorphism $\varphi $ of $\Crep$. The orbit category $C/\langle \varphi\rangle$ is isomorphic to $\Ctilde = \cdcc$. The projection $G:\Crep\to \Ctilde$ is thus 
 a Galois covering with infinite cyclic group generated by $\varphi$. It is shown in \cite{ABS3} that the corresponding pushdown functor $\textup{mod}\,\Crep \to \textup{mod}\,\Ctilde$ is always dense, so it induces an isomorphism $\zG(\textup{mod}\,\Ctilde)\cong \zG(\textup{mod}\,\Crep)/\mathbb{Z}$. Also, if $C=\End_AT$, where $T$ is a tilting module over the hereditary algebra $A$, then $\textup{mod}\,\Crep\cong\DA/\iadd(\tau_{\mathcal{D}}F^iT)_{i\in \mathbb{Z}}$, where $\tau_{\mathcal D}$ is the Auslander-Reiten translation in $\DA$ and $\iadd(\tau_{\mathcal{D}}F^iT)_{i\in \mathbb{Z}}$ is the ideal of $\DA$ consisting of all morphisms which factor through $\add(\tau_{\mathcal{D}}F^iT)_{i\in \mathbb{Z}}$. Finally, every local slice in $\zG(\textup{mod}\,\Ctilde)$ is the image under $G_\zl$ of (several) local slices in $\zG(\textup{mod}\,\Crep)$ (that is, full subquiver of $\zG(\textup{mod}\,\Crep)$ satisfying the axioms (LS1),(LS2),(LS3) of (\ref{sect 1.3}) above).  
  Throughout this paper, we identify $C_0$ with $C$, and thus any complete slice of $\textup{mod}\,C$ can be considered as a local slice in $\textup{mod}\,\Crep$.

\end{subsection}

\section{Properties of the Auslander-Reiten quiver of a cluster-tilted algebra}\label{sect 2}

\subsection{}
In this section, we let $C$ be a tilted algebra, having $\zS$ as a complete slice, and $\Ctilde=C\ltimes\DCC$ be its relation extension. The following lemma is borrowed from \cite{ADLS}; we include the proof for the convenience of the reader.

\begin{lem}\label{lem 2.1}
Let $C$ be a tilted algebra, $\zS$ a complete slice in $\textup{mod}\,C$ and $M\in\zS$, then we have:
\begin{itemize}
\item[\textup{(a)}] $M\otimes_C \DCC=0$, and
\item[\textup{(b)}] $\Hom_C(\DCC,\tau_C M) =0.$
\end{itemize}
\end{lem}
 
\begin{pf}
(a) Let $A=\End(\oplus_{X\in \zS}X)$ and $T_A$ be a tilting module such that $C=\End T_A$. Since $M\in \zS$, there exists  
 an injective $A$-module $I$ such that $M_C\cong \Hom_A(T,I).$ Using standard functorial isomorphisms, we have:
\[\begin{array}{rcl}
D(M\otimes_C\DCC) &\cong& \Hom_C(M, D\DCC)\\
 &\cong&   \Hom_C\big(\Hom_A(T,I), D\Hom_{\DA} (T,FT)     \big)           \\
  &\cong&    \Hom_C\big(\Hom_A(T,I), D\Hom_{\DA} (T,\tau^{-1}T[1])     \big)                   \\
   &\cong&       \Hom_C\big(\Hom_A(T,I), D\Hom_{\DA} (\tau T,T[1])     \big)                \\
    &\cong&              \Hom_C\big(\Hom_A(T,I), D\Ext^1_{\DA} (\tau T,T)     \big)         \\
     &\cong&         \Hom_C\big(\Hom_A(T,I), \Hom_{A} (T,\tau^2 T)     \big)              \\
      &\cong&      \Hom_A(I,t(\tau^2T)),              \\
\end{array}\]
where $t(\tau^2 T)\cong \Hom_A(T,\tau^2 T)\otimes_C T$
  is the torsion part of the $A$-module $\tau^2 T$ in the torsion pair induced by $T$ in  $\textup{mod}\,A$. Since $\tau^2 T$ is not an injective $A$-module, neither is its submodule $t(\tau^2 T)$. Since $A$ is hereditary, and $I$ is injective, we get $\Hom_A(I,\tau^2 T)=0$. 

(b) Since $\tau_CM$ precedes the complete slice $\zS$ in $\textup{mod}\,C$, it suffices to prove that $\DCC$ succedes it. Note first that
\[\begin{array}{rcl}
\DCC &\cong& \Ext^1_C(DC,\zO^{-1} C)\\ &\cong& D\sHom_C(\tau^{-1}\zO^{-1} C, DC),
\end{array}\]
using the first cosyzygy $\zO^{-1} C$ of $C$ and the Auslander-Reiten formula. Now notice that for every indecomposable summand $X$ of  $\zO^{-1} C$, there exists an injective $C$-module $J$ such that $\Hom_C(J,X)\ne 0$. But all injectives are successors of $\zS$, so there exists $L\in \zS$ such that we have a path $L\to J\to X\to*\to \tau^{-1}X$. 
This shows that every indecomposable summand  of $\tau^{-1}\zO^{-1}C$ succedes (properly) the slice $\zS$. Since no indecomposable projective module is a successor of $\zS$, we get
\[ \sHom_C(\tau^{-1}\zO^{-1} C, DC) =\Hom_C(\tau^{-1}\zO^{-1} C, DC).\]
Hence \[\DCC_C\cong D\Hom_C(\tau^{-1}\zO^{-1} C, DC)\cong \tau^{-1}\zO^{-1} C_C.\]
But as we have already shown, every indecomposable summand of $\tau^{-1}\zO^{-1} C_C$ is a (proper) successor of $\zS$. The required statement follows at once.
\qed \end{pf}

\subsection{}
\begin{prop}\label{prop 2.2}
Let $C$ be a tilted algebra, $\zS$ be a complete slice in $\textup{mod}\,C$ and $M\in \zS$. Then:
\begin{itemize}
\item[\textup{(a)}] $\tau_CM\cong\tau_{\Ctilde}M,$ and
\item[\textup{(b)}] $\tau^{-1}_CM\cong\tau_{\Ctilde}^{-1}M$.
\end{itemize}
\end{prop}

\begin{pf}
Part (a) follows directly from Lemma \ref{lem 2.1} and the main result of \cite{AZ}. Part (b) follows by duality.
\qed \end{pf}

\subsection{}\label{sect 2.3} 
We need to apply Proposition \ref{prop 2.2} also to the cluster repetitive algebra $\Crep$ of $C$.


\begin{cor}\label{cor 2.3}
Let $C$ be a tilted algebra, $\zS$ be a complete slice in $\textup{mod}\,C$ and $M\in\zS$. Then:
\begin{itemize} \item[\textup{(a)}] $\tau_C M\cong \tau_{\Crep} M$, \item[\textup{(b)}] $\tau_C^{-1} M\cong \tau_{\Crep}^{-1} M$.\qed  \end{itemize}
\end{cor}

\subsection{}
For the next lemma, we need some notations: let $A$ be a hereditary algebra, $T$ be a tilting $A$-module such that $\End_A T= C$ and $\End_{\CA} T =\Ctilde$ (where $\CA$ denotes the cluster category associated to $A$). Let also $\tilde P_x,\tilde I_x$ and $T_x$ be the indecomposable projective $\Ctilde$-module, the indecomposable injective $\Ctilde $-module and the indecomposable summand of $T$ corresponding to an object $x$ in $\Ctilde$.

\begin{lem}\label{lem 2.4} With the above notation:
\begin{itemize} \item[\textup{(a)}] For every object $x$ in $\Ctilde$, we have
$\Hom_{\CA}(T,\tau^2 T_x)\cong \tilde I_x.$
\item[\textup{(b)}] For every pair of objects $x,y$ in $\tilde C$, we have an isomorphism of the spaces of irreducible morphisms
$\textup{Irr}_{\Ctilde}(\tilde P_x,\tilde P_y)\cong \textup{Irr}_{\Ctilde}(\tilde I_x,\tilde I_y).$
\end{itemize}
\end{lem}

\begin{pf} Using standard functorial isomorphisms we have:
\[\begin{array}
{lrclll}
\textup{(a)}\qquad & \tilde I_x &\cong & D\Hom_{\CA}(T_x,T)           \\
&&\cong &   D\Hom_{\DA}(T_x,T) & \oplus&  D\Hom_{\DA}(T_x,\tau^{-1}T[1])            \\
&&\cong &    \Ext^1_{\DA}(T,\tau T_x)  &\oplus&  D\Ext^1_{\DA}(T_x,\tau^{-1}T)                    \\
&&\cong &    \Hom_{\DA}(T,\tau T[1] )  &\oplus & \Hom_{\DA}(T,\tau^{2}T_x)                    \\
&&\cong &     \Hom_{\CA}(T,\tau^{2}T_x)      . \\
\end{array}\]
\[\begin{array}{lrcl}
\textup{(b)} \qquad & \textup{Irr}_{\Ctilde} (\tilde P_x,\tilde P_y) 
&\cong& \textup{Irr}_{\CA} (T_x,T_y)\\
&&\cong& \textup{Irr}_{\CA} (\tau^2 T_x,\tau^2 T_y)\\
&&\cong& \textup{Irr}_{\Ctilde} \big(\Hom_{\CA}(T,\tau^2 T_x) ,\Hom_{\CA}(T,\tau^2 T_y)\big)\\
&&\cong&\textup{Irr}_{\Ctilde} (\tilde I_x,\tilde I_y),
\end{array}\]
where we have used the category equivalence $\Hom_{\CA}(T,-):\CA/\textup{iadd} (\tau T) \to \textup{mod}\,\Ctilde$ of \cite{BMR}, and part (a) above.
\qed
\end{pf}
\begin{rem}
Statement (b) above does not hold true in general, even for tilted algebras. Let indeed $C$ be given by the quiver
\[\xymatrix{1&2\ar[l]_\zg &3\ar[l]_\zb &4\ar[l]_\za}\] 
bound by $\za\,\zb=0.$ Note that $\textup{Irr}_C(I_1,I_2)=0$ while $\textup{Irr}_C(P_1,P_2)=k$.
\end{rem}


\section{Reflections}\label{sect 3}

\subsection{}\label{sect 3.1}
The objective of this section is to define a notion of reflection on a local slice in a cluster-tilted algebra. This will in turn induce a notion of reflection on a tilted subalgebra of the given cluster-tilted algebra.

Let, as before, $C$ be a tilted algebra, $\Ctilde=C\ltimes\DCC$ its relation-extension algebra and $\Crep$ its cluster repetitive algebra. We still identify $C$ with the full convex subcategory $C_0$ of $\Crep$. We assume throughout that $C$ is of tree type.

Let $\zG$ be a connecting component of $\textup{mod}\,C$, and $\zS$ be a complete slice in $\zG$.

Assume first that $M\in \zS$ is a source in $\zS$ which is not injective, then $(\zS\setminus\{M\})\cup\{\tau_C^{-1}M\}$ is also a complete slice in $\zG$. In the language of \cite{BOW}, these two slices are \emph{homotopic}. Homotopy is clearly an equivalence relation on slices, and there are either one or two equivalence classes in $\textup{mod}\,C$ (two if and only if $C$ is concealed). We need distinguished representatives of these classes. If there exists a complete slice in which all sources are injective $C$-modules, then such a slice is unique and is called the \emph{rightmost slice} of $\textup{mod}\,C$. We denote it as $\Sp$.
Dually, we define the \emph{leftmost slice} $\zS^-$ of $\textup{mod}\,C$. Note that, if $C$ is representation-finite, then rightmost and leftmost slices exist. 

We recall from \cite{HW} that a point $x\in C_0$ is a \emph{strong sink}  if the injective module $I_x$ has no injective module as a proper predecessor in $\textup{mod}\,C$. Clearly, strong sinks are sinks. 
The following Lemma is  obvious.

\begin{lem}\label{lem 3.1} A point  $x \in C_0$ is a strong sink if and only if  $I_x$ is an injective source of the rightmost slice $\Sp$.
\end{lem}

\begin{pf} Assume first that $I_x$ is an injective source of $\Sp$. If $x$ is not a strong sink, then there exists $y\ne x$ in $C$ such that we have a path
$I_y\leadsto I_x$. Since $\Sp$ is sincere, there exists $M\in \Sp$ and a morphism $M\to I_y$ yielding a path 
 $M\to I_y\to I_x$. Since $\Sp$  is convex in $\ind C$, we get $I_y\in \Sp$ which contradicts the hypothesis that $I_x$ is a source in $\Sp$. 
 
 Conversely, assume $x$ to be a strong sink in $C$, and suppose that $I_x$
is not an injective source of $\Sp$. Because $\Sp$ is sincere, then there exist $N\in \Sp$ and a morphism $N\to I_x$. Now there exists a source (necessarily injective) $I_z$ in $\Sp$ and a path $I_z\leadsto N$ in $\Sp$. This yields a path $I_z\leadsto N\to I_x$, contrary to the hypothesis. 
 \qed
\end{pf}
 
\subsection{The completion $G_x$}\label{sect 3.2}
Let $x$ be a strong sink in $C$.
We define the \emph{completion} $G_x$ of $x$ in $\Sp$ to be a non-empty
full connected subquiver of $\Sp$ such that
\begin{itemize} \item[\textup{(a)}] $I_x\in G_x$,
\item[\textup{(b)}] $G_x$ is closed under predecessors in \Sp,
\item[\textup{(c)}] If $I\to M$ is an arrow in \Sp, with $I\in G_x$ injective, then $M\in G_x$,
\item[\textup{(d)}] If $N\to I$ is an arrow in \Sp, with $I\in G_x$ injective, then $N $ is injective (and in $G_x$).
 \end{itemize}

Completions do not always exist. 
\begin{exmp}
The tilted algebra $C$ given by the quiver
\[\xymatrix@C40pt{1&2\ar@<2pt>[l]^\zg\ar@<-2pt>[l]_\zb  &3\ar[l]_\za}\]
bound by $\za\,\zb=0$ admits the complete rightmost slice consisting of the modules $I_1, S_2$ and $I_2$, and $I_1$ is the only source. A part of the Auslander-Reiten quiver of $\textup{mod}\,C$ containing this slice is shown below, where modules are represented by their dimension vectors.
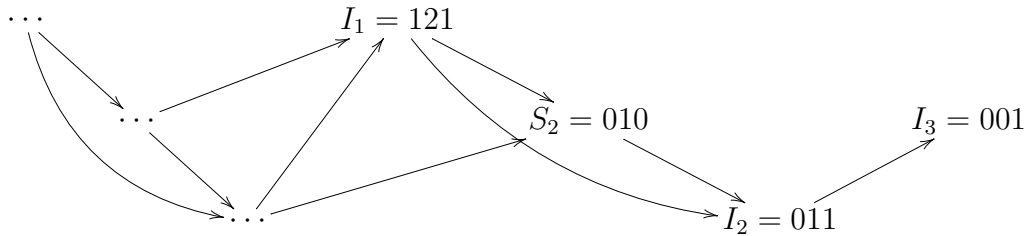
\[\xymatrix{\cdots\ar[rd]\ar@/_20pt/[rrdd] &&& I_1=121 \ar[rd]\ar@/_20pt/[rrdd] & & \\
&\cdots\ar[rd]\ar[rru] &&& S_2 =010\ar[rd] && I_3=001\\
&&\cdots\ar[rru]\ar[ruu]&&&I_2=011\ar[ru]}
\]
In this example $G_1$ does not exist, because by condition (c) it would contain both $S_2 $ and $I_2$, and this contradicts condition (d).
\end{exmp}
The tilted algebra $C$ in the example above is of euclidean type $\tilde{\mathbb{A}}_2$, so it is not of tree type. The following Lemma guarantees the existence of some completion in a rightmost slice, if the tilted algebra is of tree type.

\begin{lem} \label{lem 3.2}
Let $C$ be a tilted algebra of tree type having a rightmost slice $\Sp$. Then there exists a  strong sink $x$ in $C$ such that the completion $G_x$ exists.
\end{lem}
\begin{pf}
Let $I_{x_1}$ be a source in $\Sp$ and $G_1'$ its closure under condition (c) above, then let $G_1$ be the closure of $G_1'$ under condition (b).

If $G_1$ satisfies condition (d), then we are done. Otherwise there exist an injective $I\in G_1$ and an arrow $N\to I$ in  $\Sp$ with $N$ not injective. Then there exists a sectional path in $\Sp$ ending at $N$. Let $I_{x_2}$ be the source of such a path.

Let $G_2'$ be the closure of $I_{x_2}$ under condition (c), and then
let $G_2$ be the closure of $G_2'$ under condition (b). Clearly, $G_2'$  does not contain the injective $I$, since there is an arrow $N\to I $ in the sectional path, with $N$ non-injective. Using that $\Sp$ is a tree, we see that $I_{x_1}\notin G_2$.

If $G_2$ satisfies condition (d), then we are done. Otherwise we repeat the procedure. Since $\Sp$ is a tree, this procedure must ultimately stop.
\qed\end{pf}

\begin{exmp}
Let $C$ be the tilted algebra of tree type $\mathbb{D}_5$ given by the quiver
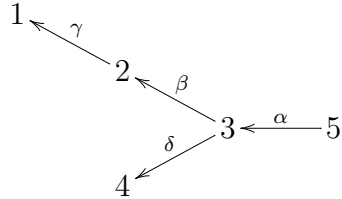
\[\xymatrix@R10pt@C30pt{1\\&2\ar[ul]_\zg \\&&3\ar[ul]_\zb\ar[dl]_\zd & 5\ar[l]_\za \\&4
}
\]
bound by $\za\,\zb\,\zg=0$ and $\za\,\zd =0$.
Its Auslander-Reiten quiver is shown below.
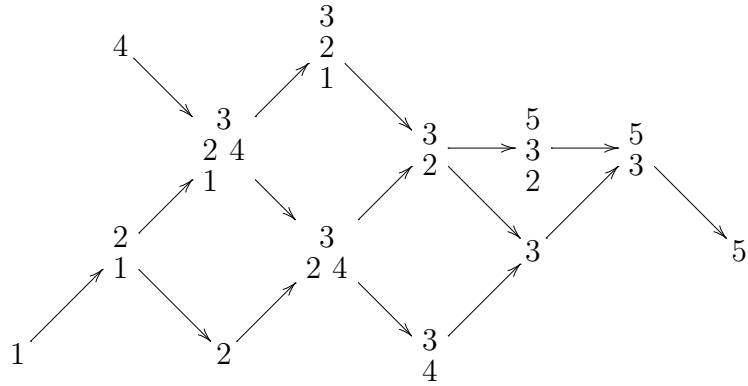
\[\xymatrix@!R=15pt@!C=15pt{&4\ar[rd]
&&{\begin{array}{c} 3  \vspace{-10pt}\\ 2 \vspace{-10pt}\\ 1\end{array}}\ar[rd]
\\
&&{\begin{array}{c}3\vspace{-10pt}\\ 2\ 4  \vspace{-10pt}\\ 1\hspace{10pt}  \end{array}}\ar[rd]\ar[ru]
&&{\begin{array}{c} 3  \vspace{-10pt}\\ 2\end{array}}\ar[r]\ar[rd]
&{\begin{array}{c} 5  \vspace{-10pt}\\ 3 \vspace{-10pt}\\ 2\end{array}}\ar[r]
&{\begin{array}{c} 5  \vspace{-10pt}\\ 3 \end{array}}\ar[rd]
\\
&{\begin{array}{c} 2  \vspace{-10pt}\\ 1\end{array}}\ar[rd]\ar[ru]
&&{\begin{array}{c} 3  \vspace{-10pt}\\ 2\ 4\end{array}}\ar[rd]\ar[ru]
&&3\ar[ru] 
&&5
\\
1\ar[ru]&& 2\ar[ru]
&&{\begin{array}{c} 3  \vspace{-10pt}\\ 4\end{array}}\ar[ru]
}
\] (here, modules are represented by their composition factors).
The rightmost slice 
\[\left\{ {\begin{array}{c} 3  \vspace{-10pt}\\ 2 \vspace{-10pt}\\ 1\end{array}}\,,\,
{\begin{array}{c} 3  \vspace{-10pt}\\ 2\end{array}}\,,\,
{\begin{array}{c} 3  \vspace{-10pt}\\ 4\end{array}}\,,\,
{\begin{array}{c} 5  \vspace{-10pt}\\ 3 \vspace{-10pt}\\ 2\end{array}}\,,\,
3
\right\}\]
in this example has the two injective sources: $I_1$ and $I_4$. We have 
\[G_1=\left\{ {\begin{array}{c} 3  \vspace{-10pt}\\ 2 \vspace{-10pt}\\ 1\end{array}}\,,\,
{\begin{array}{c} 3  \vspace{-10pt}\\ 2\end{array}}\right\} \quad \textup{ and } \quad G_4 =  
\left\{ {\begin{array}{c} 3  \vspace{-10pt}\\ 2 \vspace{-10pt}\\ 1\end{array}}\,,\,
{\begin{array}{c} 3  \vspace{-10pt}\\ 2\end{array}}\,,\,
{\begin{array}{c} 3  \vspace{-10pt}\\ 4\end{array}}\,,\,
3
\right\}.\]

\end{exmp}
\subsection{The reflection of a slice}\label{sect 3.3}
Let now $x $ be a sink in $C$ such that the completion $G_x$ exists. We   then say that $x$ is an \emph{admissible sink}. We are now able to define the reflection $\zS'=\sx\Sp$ of the complete slice $\Sp$. The set of objects in $G_x$ is of the form $\mathcal{J}\sqcup\mathcal{M}$, where $\mathcal{J}$ and $\mathcal {M}$ consist respectively of the injective, and the non-injectives in $G_x$. Let $\mathcal{P}=\{P_x\in \textup{mod}\,C_1\mid I_x\in \mathcal{J}\}$, where we recall that $C_1$ is the copy of $C$ next to $C_0$ on the diagonal blocks of $\Crep$. We then set
\[\sx\Sp= \big(\Sp\setminus G_x\big) \cup \mathcal{P}\cup\tau_{\Crep}^{-1}\mathcal{M}.\]
Recall that, by Corollary \ref{cor 2.3}, $\tau_{\Crep}^{-1} M \cong \tau_{C}^{-1} M$ for every $M\in \Sp$.

\begin{lem}\label{lem 3.3}
$\sx\Sp$ is a local slice in $\textup{mod}\,\Crep$.
\end{lem}

\begin{pf}
We first consider in the cluster category $\CA$ the full subquiver defined by:
\[\zS''=\big(\Sp\setminus G_x\big)\cup\tau_{C}^{-1}\mathcal{M} \cup\tau^{-1}_{\CA} \mathcal{I}.
\]
Thus $\zS''$ is a local slice in $\CA$ because $G_x$ is closed under predecessors and we have $\zS'=(\zS''\setminus  \tau^{-1}_{\CA} \mathcal{I}) \sqcup \mathcal{P}.$

We   claim that $\zS'$ is connected. The objects lying in $\zS'$ and $\zS''$ are in one-to-one correspondence, since any object of $\zS'$ is either an object of $\zS''$ or the Auslander-Reiten translate of an object in $\zS''$. Hence it is enough to show that whenever there is an arrow between $M'',N''$ in $\zS''$, then there is an arrow between the two corresponding objects $M',N'$ in $\zS'$.

Because of Lemma \ref{lem 2.4}\,(b), we only need to consider the case where $M''\in(\Sp\setminus G_x)\cup \tau^{-1}_{C}\mathcal{M} $ and $N''\in \tau^{-1}_{\CA}\mathcal{I}$. Thus $M'=M'' $ and $N'=\tau^{-1}_{\CA} N''=\tau^{-2}_{\CA} I$ for some $I\in\mathcal{I}\subset G_x$.

 Either we have $M''\to N''$ or $N''\to M''$ in $\zS''$. In the latter case, there is an arrow from $(M'=M'')$ to $(N'=\tau^{-1}_{\CA}N'')$ in $\zS'$, and we are done. On the other hand, if $M''\to N''$,
 then there is an arrow $\tau_{\CA} N'' \to M''$ with $\tau_{\CA}N''=I\in G_x$ injective, and thus $M'\in G_x$, by condition (c) of the completion $G_x$. This establishes our claim.

Consequently, $\zS'$ may be identified to a local slice in $\DC$. Since, furthermore, $\zS'$ consists of $\Crep$-modules then, by \cite{ABS3}, $\zs'$ is a local slice in $\textup{mod}\,\Crep$.\qed
\end{pf}

\subsection{A hereditary subcategory}
 We   deduce from our definition of reflection of $\Sp$ a definition of reflection of the tilted algebra $C$, which we   denote by $\sC$.
 
 Define $\mathcal{S}_x$ to be the full subcategory of $C$ consisting of the objects $y$ such that $I_y\in G_x$.
 
\begin{lem}
\label{lem 3.4}
With the above notation
\begin{itemize}
 \item[\textup{(a)}]  $\mathcal{S}_x$ is hereditary,
\item[\textup{(b)}] $\mathcal{S}_x$ is closed under successors in  $C$,
\item[\textup{(c)}] $C$ may be written in the form 
\[C= \left[ \begin{array}{cc} \ H\ &\ 0\ \\M&C'
\end{array}\right]\]
with $H$ hereditary, $C'$ tilted and $M$ a $C'$-$H$-bimodule.
\end{itemize}
\end{lem}
\begin{pf}
(a) We let  $H=\End(\oplus_{y\in\mathcal{S}_x} I_y)$. Then $H$ is a full subcategory of the hereditary algebra $A=\End (\oplus_{X\in\Sp} X)$. Therefore $H$ is also hereditary, that is, $\mathcal{S}_x$ is hereditary.

(b) Let $y\in \mathcal{S}_x$ and $y\to z$ be an arrow in $C$. Then there exists a morphism $I_z\to I_y$. Since $I_z$ is an injective $C$-module and $\Sp$ is sincere, there exist $N\in \Sp$ and a morphism $N\to I_z$. Thus we have $N\to I_z\to I_y$. Since $N,I_y\in \Sp$ and $\Sp$ is convex in $\textup{mod}\,C$, then $I_z\in \Sp$ and so $z\in \mathcal{S}_x$.

(c) This follows at once from (a) and (b).\qed 
\end{pf}
\subsection{The structure of the cluster duplicated algebra} 
We recall from \cite{ABS3} that the cluster duplicated algebra $\Cbar$ of $C$ is the (finite dimensional) matrix algebra
\[\Cbar   \ = \ \left[ \begin{array}{cc} C&0\\ \ \DCC\ &\ C\ 
\end{array}
\right]
\]
with the ordinary matrix addition and the multiplication induced from that of $C$ and from the $C$-$C$-bimodule structure of $\DCC$. Clearly, $\Cbar$ is useful as a ``building block'' for the cluster repetitive algebra $\Crep$.

\begin{cor}
\label{cor 3.5}
The cluster duplicated algebra of $C$ is of the form 
\[\Cbar\ =\ \left[ \begin{array}{cccc}\ \  H\ \ &\ \ 0\ \  &\ \ 0\ \ &\ \ 0\ \ \\ M&C'&0&0\\ 0&F_0&H&0 \\0&F_1&M&C'
\end{array}
\right],
\]
where $F_0= \Ext^2_C(DC',H)$ and $F_1=\Ext^2_C(DC',C').$
\end{cor}
\begin{pf}
 We start by writing $C$ in the matrix form of Lemma \ref{lem 3.4}\,(c). Since, by definition, $H$ consists of the objects $y$ in $C$ such that $I_y\in G_x\subset \Sp$, then the projective dimension $\pd_C DH$ is at most $1$, hence $\Ext^2_C(DH,-)=0$. The result follows upon multiplying by idempotents. \qed 
\end{pf}

\subsection{The reflection of a tilted algebra} 
We can now define the \emph{reflection} $\sC$ of $C$ to be the matrix algebra
\[ \sC =
\left[ \begin{array}{cc} \ C'\ & \ 0\ \\ F_0&H
\end{array}
\right],\]
where $F_0= \Ext^2_C(DC',H)$. Note that $\sC$ is a quotient algebra of $\Crep$.

We now prove that this definition is compatible with the definition of reflection of local slices. We recall that the \emph{support} $\Supp \mathcal{X}$ of a subclass $\mathcal{X}$ of $ \Crep$ is the full subcategory of $\Crep$ having as objects the $x$ in $\Crep$ such that there exists a module $M\in\mathcal{X}$ satisfying $M(x)\ne 0$.

\begin{prop}
\label{prop 3.6} The reflection
$\sC$ is a tilted algebra having $\sx\Sp$ as a complete slice. Moreover, the cluster-tilted algebras of $C$ and $\sC$ and the cluster repetitive algebras of $C$ and $\sC$ are isomorphic. 
%
\end{prop}

\begin{pf}
It follows directly from the definition of $\sx\Sp$ that $\Supp(\sx\Sp)\subset \sC$. Indeed, in the notation of Lemma \ref{lem 3.3}, we have $\sx\Sp=(\Sp\setminus G_x)\cup\mathcal{P}\cup\tau^{-1}_{\Crep}\mathcal{M}.$ Since, as observed before, $\tau^{-1}_{\Crep}\mathcal{M}
\cong\tau^{-1}_{C}\mathcal{M}$ by Corollary \ref{cor 2.3}, and the injectives in $\mathcal{I}$ are replaced by the projectives in $\mathcal{P}$, then we get the wanted inclusion.

Now, as shown in Lemma \ref{lem 3.3}, $\sx\Sp$ is a local slice in $\textup{mod}\,\Crep$. Denoting by $G_\zl:\textup{mod}\,\Crep \to \textup{mod}\,\Ctilde$ the pushdown functor associated to the Galois covering $G:\Crep\to \Ctilde$, we get that $G_\zl (\sx\Sp)$ is a local slice in $\textup{mod}\,\Ctilde$. By \cite{ABS2}, $C^*=\Ctilde /\Ann(G_\zl(\sx\Sp))$ is a tilted algebra of the same type as $C$. Moreover we have $\Ctilde= C\ltimes \DCC \cong C^*\ltimes\Ext^2_{C^*}(DC^{*},C^*)$ so that we also have $\Crep=\stackrel{\scriptscriptstyle\vee}{C^*}$.

On the other hand, $\sx\Sp$ is a complete slice in $\textup{mod}\,C^*$ so, in particular, it is sincere over $C^*$. Therefore, $\Supp \sx\Sp=C^*$. Using that $\Crep=\stackrel{\scriptscriptstyle\vee}{C^*}$, we thus have $C^*\subset \sC$. Finally, since the  Grothendieck groups of $C^*, \sC$ and $C$ are all of the same rank, it follows that the full subcategories $C^*$ and $\sC$ of $\Ctilde $ are equal. This completes the proof. \qed
\end{pf}

Dually, one defines coreflections $\zs_x^-$ with respect to admissible sources $x$. We leave the straightforward statements to the reader.


\section{Main result}\label{sect 4}
\subsection{The distance between two local slices} 
We introduce the following notation. Let  $\zS_1,\zS_2$ be two local slices in $\textup{mod}\,\Crep$, considered as embedded in $\DC$. We define $\dc(\zS_1,\zS_2)$ to be the number of $\tau F^j T_i$ (where $1\le i\le \textup{rk} K_0(C)$ and $j\in \mathbb{Z}$) in $\DC$ such that either $\zS_1< \tau F^j T_i < \zS_2$, or  $\zS_2< \tau F^j T_i < \zS_1$. 

Note that  $\dc(\zS_1,\zS_2)$ is always a non-negative integer but it can be arbitrarily large. Also, if $\Crep$ is locally representation-finite (that is, $\Ctilde$ is representation-finite), then  $\dc(\zS_1,\zS_2)=0$ if and only if the local slices $G_\zl\zS_1$ and $G_\zl \zS_2$ in $\textup{mod}\,\Ctilde $ are homotopic in the sense of \cite{BOW} (see section (\ref{sect 3.1}) above).

\begin{lem}
\label{lem 4.1}
Let $\zS_1,\zS_2,\zS_3$ be local slices in $\textup{mod}\,\Ctilde$, then:
\begin{itemize} \item[\textup{(a)}]  $\dc(\zS_1,\zS_2)=\dc(\zS_2,\zS_1),$
\item[\textup{(b)}] $\dc(\zS_1,\zS_3)\le \dc(\zS_1,\zS_2)+\dc(\zS_2,\zS_3)$.
\end{itemize}
\end{lem} 
 
\begin{pf}
(a) is obvious and (b) follows from a straightforward counting argument. \qed
\end{pf}

\subsection{The metric space of fibre quotients of a cluster repetitive algebra}
Clearly, $\dc$ is not yet a distance function. Our objective is to use it in order to define a distance function. We say that an algebra $C'$ is a \emph{fibre quotient} of $\Crep$ if $C'$ is tilted and such that $\stackrel{\scriptscriptstyle\vee}{C'}\cong \Crep$. This terminology is motivated by the observation  that such an algebra $C'$ lies in the fibre of $\Crep$ under the mapping $C\mapsto \Crep$ from the class of tilted algebras to the class of cluster repetitive algebras.

Let now $C_1,C_2$ be two fibre quotients of $\Crep$, and $\zS_1,\zS_2 $ be complete slices in $\textup{mod}\,C_1,\textup{mod}\,C_2$ respectively, considered as local slices in $\textup{mod}\,\Crep$. Then we set
\[\dc(C_1,C_2)=\dc(\zS_1,\zS_2).\]
This does not depend on the choice of the complete slices $\zS_1$ and $\zS_2$. Indeed, let $\zS_1,\zS_1'$ be two complete slices in $\textup{mod}\, C_1$, then it is clear that $\dc(\zS_1,\zS_1')=0$. Hence Lemma \ref{lem 4.1}\,(b) yields 
$\dc(\zS_1,\zS_2)\le \dc(\zS_1,\zS_1')+\dc(\zS_1',\zS_2) =\dc(\zS_1',\zS_2)$. Similarly, 
$\dc(\zS_1',\zS_2)\le \dc(\zS_1,\zS_2)$, so $\dc(\zS_1,\zS_2)=\dc(\zS_1',\zS_2)$, and our notion is well-defined.

\begin{prop}
\label{prop 4.2} Let $C_1,C_2$ be two fibre quotients of $\Crep$, then $\dc(C_1,C_2)=0$ if and only if $C_1=C_2$.
\end{prop}

\begin{pf}
Assume indeed that $\dc(C_1,C_2)=0. $ Let $\zS_1,\zS_2$ be complete slices in $\textup{mod}\,C_1, \textup{mod}\,C_2$, respectively, considered as local slices in $\textup{mod}\,\Crep$. By \cite{ABS2}, we have $C_1=\Crep/\Ann \zS_1$ and  $C_2=\Crep/\Ann \zS_2$.

Let $T$ be a tilting module over the hereditary algebra $A$ such that $\End_A T\cong C$, and $\End_{\CA} T\cong \Ctilde$, (so that $\End_{\DA}(\oplus_{i\in \mathbb{Z}} F^i T)=\Crep$).
By \cite{ABS2}, the annihilator $\Ann \zS_1$ is generated by the arrows $\za:(x_0,i)\to(y_0,j)$ of $\Crep$ (here $x_0,y_0$ are points of $C_1$, while $i,j\in \mathbb{Z}$) such that the corresponding morphism $f_\za:F^j T_{y_0} \to F^i T_{x_0}$ in the derived category lies in $\Hom_{\DA}(F^j T, F^{j+1} T)$ and $\zS_1=F^j DA$. 
Now, this is the case if and only if
\[ F^jT_{y_0}\le \zS_1\le \tau^2F^{j+1} T_{x_0}\]
in $\DA$. Indeed, notice first that the existence of the arrow $\za$ means that $i\in\{j,j+1\}$. Moreover $\tau^2FT_{x_0}=\tau T_{x_0}[1]\ge DA$ implies $\tau^2F^{j+1}T_{x_0}\ge F^j DA =\zS_1$. On the other hand, $T_{y_0}\le DA$ gives clearly $F^jT_{y_0} \le F^j DA=\zS_1$.

We next claim that $\dc(\zS_1,\zS_2)=0 $ implies 
\[F^jT_{y_0} \le \zS_2\le \tau^2F^{j+1}T_{x_0}.\]
Indeed, if $F^jT_{y_0} \nleq\zS_2$, then $\zS_2 <F^jT_{y_0}$, so that  $\zS_2 <\tau F^jT_{y_0}$ because $\tau F^jT_{y_0}\notin \zS_2$. This implies that $\zS_2 <\tau F^jT_{y_0}< \zS_1$ and we have a contradiction to $\dc(\zS_1,\zS_2)=\dc(C_1,C_2)=0$. 
On the other hand, if $ \zS_2\nleq \tau^2 F^{j+1} T_{x_0}$, then  $\tau^2 F^{j+1} T_{x_0}  < \zS_2$ and so  $\tau F^{j+1} T_{x_0}  < \zS_2$ because $ \tau F^{j+1} T_{x_0}  \notin \zS_2$. This implies that $\zS_1<  \tau F^{j+1} T_{x_0}  < \zS_2$, another contradiction to $\dc(\zS_1,\zS_2)=\dc(C_1,C_2)=0$. 
This establishes our claim.

Now, that claim implies that the annihilators of $\zS_1$ and $\zS_2$ have the same generators. Therefore $C_1=C_2$. Since the converse is obvious, the proof of the proposition is complete.\qed
\end{pf}

\begin{cor}
\label{cor 4.3} The set $\stackrel{\scriptscriptstyle\vee}{\mathcal{F}}$ of all fibre quotients of $\Crep$ is a discrete metric space with the distance $\dc$.
\end{cor}
\begin{pf}
 It follows from Lemma \ref{lem 4.1} and Proposition \ref{prop 4.2} that $\dc$ is a distance in $\stackrel{\scriptscriptstyle\vee}{\mathcal{F}}$. It is clear that the resulting metric space is discrete.\qed
\end{pf}
\subsection{The metric space of fibre quotients of a cluster-tilted algebra}
We now bring down this information to $\Ctilde$. We say that an algebra $C'$ is a \emph{fibre quotient} of $\Ctilde$ if $C'$ is tilted and such that $\tilde{C'}\cong \Ctilde$. 
Let $C_1,C_2$ be two fibre quotients of $\Ctilde$, then we set
\[d (C_1,C_2)= \operatorname*{min}_{C_1^*,C_2^*\in\stackrel{\scriptscriptstyle\vee}{\mathcal{F}}}\ \{\dc(C_1^*,C_2^*)\mid GC_1^*=C_1, GC_2^*=C_2\}.\]

\begin{cor}
\label{cor 4.4} Let $C_1,C_2$ be two fibre quotients of $\Ctilde$, then $d(C_1,C_2)=0$ if and only if $C_1=C_2$.
\end{cor}
\begin{pf}
This follows immediately from Proposition \ref{prop 4.2}.\qed
\end{pf}
\begin{rem}
This gives another interpretation and proof of \cite[Theorem 4.13]{BOW}.
\end{rem}

Notice that while our definition implies that the set $   \stackrel{\scriptscriptstyle\vee}{\mathcal{F}}$ of fibre quotients of $\Crep$ is infinite, clearly the set $\tilde{\mathcal{F}}$ of fibre quotients of $\Ctilde$ is finite. Moreover, it is easily seen that $   \stackrel{\scriptscriptstyle\vee}{\mathcal{F}}$ is (trivially) a topological covering of $\tilde{\mathcal{F}}$.
\begin{cor}\label{cor 4.6}
The set $\tilde{\mathcal{F}}$ of all fibre quotients of $\Ctilde$ is a discrete metric space with the distance $d$.
\end{cor}
\begin{pf}
This follows from Corollary \ref{cor 4.3}.\qed
\end{pf}

\subsection{}
The following lemma and its proof, which relate fibre quotients of $\Ctilde$ and $\Crep$, are valid without assuming that $C$ is of tree type.
\begin{lem}\label{lem 4.5}
Let $C$ be a tilted algebra. If $C'$ is a fibre quotient of $\Ctilde$, then $G^{-1}(C)$ is the $\varphi$-orbit of a fibre quotient of $\Crep$.
Conversely, if $C^*$ is a fibre quotient of $\Crep$, then $G(C^*)$ is a fibre quotient of $\Ctilde$.
\end{lem}
\begin{rem}
By abuse of language, we quote from now on this lemma by saying that $C'$ is a fibre quotient of $\Ctilde$ if and only if $C'$ is a fibre quotient of $\Crep$.
\end{rem}
\begin{pf}
Suppose $\Crep=\ \stackrel{\scriptscriptstyle\vee}{C^*}$. Let $\zS$ be a complete slice in $\textup{mod}\,C$ considered as a local slice in $\Crep=\ \stackrel{\scriptscriptstyle\vee}{C^*}$. By \cite{ABS2}, $\zS$ lifts isomorphically as a section both in $\DC$ and in $\mathcal{D}^b(\textup{mod}\, C^*)$. This implies that we have equivalences of triangulated categories $\phi : \DC \stackrel{\cong}{\to}\DkS$ and $\phi^* :\mathcal{D}^b(\textup{mod}\, C^*) \stackrel{\cong}{\to}\DkS$. Let $T=\phi\, C$ and $T^*=\phi^*C^*$. Then:
 \[\begin{array}{rcl}
\End_{\DkS} (\oplus_{j\in \mathbb{Z}} F^j T) &\cong&  \End_{\DC} (\oplus_{j\in \mathbb{Z}} F^j C) \\
&\cong&   \Crep           \\ 
&\cong&    \stackrel{\scriptscriptstyle\vee}{C^*}          \\ 
&\cong&      \End_{\DCprime} (\oplus_{j\in \mathbb{Z}} F^j C^*)        \\ 
&\cong&        \End_{\DkS} (\oplus_{j\in \mathbb{Z}} F^j T^*)  .    
\end{array}
 \]
Define $C'=G(C^*)$, then, passing to the cluster category, we have 
 $\CC\cong \mathcal{C}_{k\zS} \cong \mathcal{C}_{C'}$ 
 and
 \[\begin{array}
{rcl}
\Ctilde &\cong & \End_{\CC} C \\
&\cong & \End_{\mathcal{C}_{k\zS}} T \\
&\cong & \End_{\mathcal{C}_{k\zS}} T^* \\
&\cong & \End_{\mathcal{C}_{C'}} C' \\
&\cong& \tilde{C'}.
\end{array}
 \]
 This proves the sufficiency. The necessity is obvious.\qed
\end{pf}

\subsection{Example}\label{ex 4.7}
Let $\Ctilde$ be the cluster-tilted algebra of type $\mathbb{A}_5$ given by the quiver 

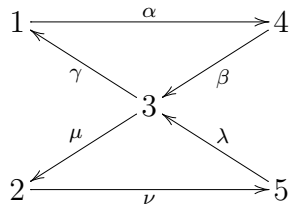
\[
\xymatrix@R=20pt@C=40pt{
1\ar[rr]^\za&&4\ar[ld]^\zb\\
&3\ar[lu]^\zg\ar[ld]_\mu\\
2\ar[rr]_\nu&&5\ar[lu]_\zl
}
\]
bound by  $\za\,\zb=0$, $\zb\,\zg=0$, $\zg\,\za=0$ $\zl\,\mu=0$, $\mu\,\nu=0$ and $\nu\,\zl=0$. 
Its Auslander-Reiten quiver is shown in
Figure \ref{fig arq}, where modules are represented by their Loewy series and we identify the vertices that have the same label, thus creating a Moebius strip.
\begin{figure}
\[
\xymatrix@!R=15pt@!C=15pt{
 {\begin{array}{c}4 \vspace{-10pt} \\ \vspace{-10pt} 3\\2\end{array}} \ar[rd]
 &&&& {\begin{array}{c}2 \vspace{-10pt} \\ 5\end{array}} \ar[rd]
 &&&& {\begin{array}{c}5 \vspace{-10pt} \\ \vspace{-10pt} 3\\1\end{array}}
 \\
&{\begin{array}{c}4 \vspace{-10pt} \\3\end{array}}\ar[rd]
&&\ 5\ \ar[ru]&& \ {2}\ar[rd]
&&{\begin{array}{c}3 \vspace{-10pt} \\1\end{array}}\ar[ru]\ar[rd]
\\
 {\begin{array}{c}3\end{array}}\ar[ru]\ar[rd]
&&{\begin{array}{c}4\ 5 \vspace{-10pt} \\3\end{array}}\ar[ru]\ar[rd]
&&&&{\begin{array}{c}3 \vspace{-10pt} \\1\ 2\end{array}}\ar[ru]\ar[rd]
&&{\begin{array}{c}3\end{array}}
\\
&{\begin{array}{c}5 \vspace{-10pt} \\3\end{array}}\ar[ru]&&4 \ar[rd]
&&1\ar[ru]
&&{\begin{array}{c}3 \vspace{-10pt} \\2\end{array}}\ar[ru]\ar[rd]
\\
{\begin{array}{c}5 \vspace{-10pt} \\3 \vspace{-10pt} \\1\end{array}}\ar[ru]
&&&&{\begin{array}{c}1 \vspace{-10pt} \\4\end{array}}\ar[ru]
&&&&{\begin{array}{c}4 \vspace{-10pt} \\3 \vspace{-10pt} \\2\end{array}}
}
\]
\caption{Auslander-Reiten quiver of Example \ref{ex 4.7}}\label{fig
  arq}
\end{figure}
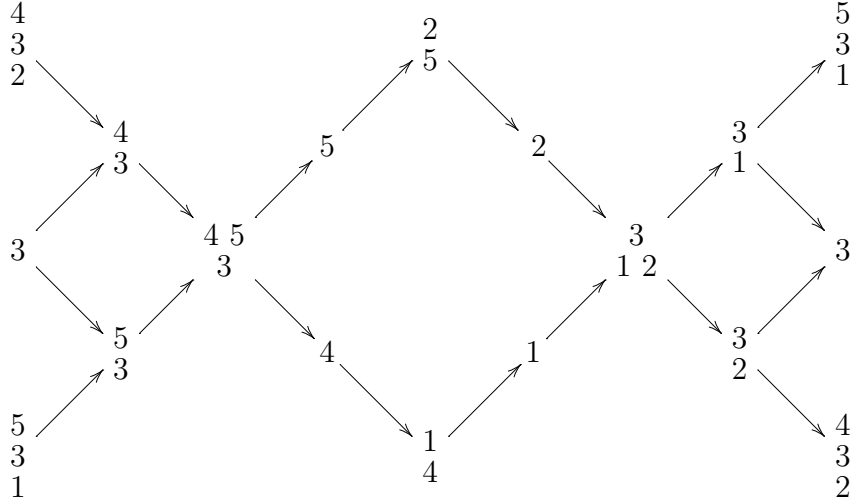
Let $\zS_1,\zS_2,\zS_3 $ be respectively given by 
\[\begin{array}{rcl}
\zS_1&=&\left\{ \ {\begin{array}{c}4 \vspace{-10pt} \\ \vspace{-10pt} 3\\2\end{array}}\ ,\ 
{\begin{array}{c}4 \vspace{-10pt} \\3\end{array}}\ ,\ 
{\begin{array}{c}4\ 5 \vspace{-10pt} \\3\end{array}}\ ,\  
4\ ,\ 
{\begin{array}{c}1 \vspace{-10pt} \\4\end{array}}
\ \right\} 
\\
\zS_2&=&\left\{\ {\begin{array}{c}5 \vspace{-10pt} \\3 \vspace{-10pt} \\1\end{array}}\ ,\ 
{\begin{array}{c}5 \vspace{-10pt} \\3\end{array}}\ ,\ 
{\begin{array}{c}4\ 5 \vspace{-10pt} \\3\end{array}}\ ,\  
5\ ,\ 
{\begin{array}{c}2 \vspace{-10pt} \\ 5\end{array}}
\ \right\}
\\
\zS_3 &=&\left\{\ {\begin{array}{c}2 \vspace{-10pt} \\ 5\end{array}}\ ,\ 
2\ ,\ 
{\begin{array}{c}3 \vspace{-10pt} \\1\ 2\end{array}}\ ,\ 
{\begin{array}{c}3 \vspace{-10pt} \\2\end{array}}\ ,\ 
{\begin{array}{c}4 \vspace{-10pt} \\3 \vspace{-10pt} \\2\end{array}}
\ \right\}.
\end{array}
\]
Then $C_1=\Ctilde/\Ann \zS_1 $ is given by the quiver
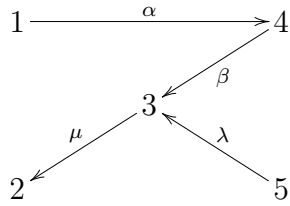
\[\xymatrix@R=20pt@C=40pt{
1\ar[rr]^\za&&4\ar[ld]^\zb\\
&3\ar[ld]_\mu\\
2&&5\ar[lu]_\zl
}
\]
while $C_2=\Ctilde/\Ann \zS_2 $ is given by the quiver 
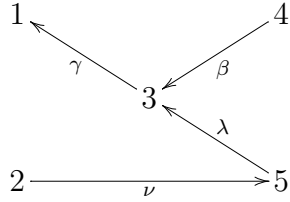
\[
\xymatrix@R=20pt@C=40pt{
1&&4\ar[ld]^\zb\\
&3\ar[lu]^\zg\\
2\ar[rr]_\nu&&5\ar[lu]_\zl
}
\]
and $C_3=\Ctilde/\Ann \zS_3 $ is given by the quiver 
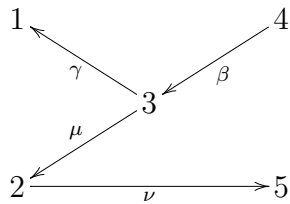
\[
\xymatrix@R=20pt@C=40pt{
1&&4\ar[ld]^\zb\\
&3\ar[lu]^\zg\ar[ld]_\mu\\
2\ar[rr]_\nu&&5
}
\]
with the inherited relations in each case. Then we have $d(C_1,C_2)=d(C_1,C_3)= d(C_2,C_3)=2.$
Notice that if $\Ctilde$ has $n$ points, then clearly, for any two fibre quotients $C_1,C_2$ of $\Ctilde$, we have $d(C_1,C_2)\le \lfloor \frac{n}{2}\rfloor$.

\subsection{}
 We are now able to state and prove the key lemma.
\begin{lem}\label{lem 4.8}
Let $\zS_1,\zS_2 $ be two local slices in the same transjective component of $\textup{mod}\,\Crep$ such that $\dc(\zS_1,\zS_2)\ne 0$. Then either:
\begin{itemize} \item[\textup{(a)}] there exists a rightmost slice $\zS_1^+$ such that $\dc(\zS_1,\zS_1^+)=0$ and a reflection $\sx$ such that $\dc(\sx\zS_1^+,\zS_2)<\dc(\zS_1,\zS_2),$ or 
\item[\textup{(b)}] there exists a leftmost slice $\zS_1^-$ such that $\dc(\zS_1,\zS_1^-)=0$ and a coreflection $\zs_y^-$ such that $\dc(\zs_y^-\zS_1^-,\zS_2)<\dc(\zS_1,\zS_2).$ 
 \end{itemize}
\end{lem}
\begin{pf}
(1) Assume first that $\zS_1\cap \zS_2 =\emptyset$, then we can assume without loss of generality that $\zS_1<\zS_2$. Let $\zS_1^+$ be the rightmost slice such that $\dc(\zS_1,\zS_1^+)=0$. Such a rightmost slice exists since $\dc(\zS_1,\zS_2)\ne 0$ and the two slices lie in the same transjective component. Let $x=(x_0,j)$ be an admissible sink in $\zS_1^+$. We claim that $\sx\zS_1^+$ gives the result. Indeed, $T_{x_0}$ is such that 
\[ \zS_1 < \tau F^j T_{x_0} <\zS_2\]
in $\DC$, but $\tau F^j T_{x_0} <\sx\zS_1^+$. Also, if $T_{y_0}$ is such that $\sx\zS_1^+<\tau F^i T_{y _0} < \zS_2  $ in \DC, then $\zS_1\le\zS_1^+<\tau F^i T_{y _0} < \zS_2$. Moreover, $\zS_2 <\tau F^i T_{y _0} <\sx\zS_1^+$ is impossible, because $\sx\zS_1^+\le \zS_2$. We deduce that $d(\sx\zS_1^+,\zS_2)<d(\zS_1,\zS_2)$.
This proves (a). Similarly, assuming $\zS_2<\zS_1$ yields (b).

(2)  Suppose now that $\zS_1\cap \zS_2\ne\emptyset$. Since $\dc(\zS_1,\zS_2)\ne 0$, there exists $z=(z_0,j)$ such that either $\zS_1 <\tau F^j T_{z_0} <\zS_2 $ or  $\zS_2 <\tau F^j T_{z_0} <\zS_1 $. Assume  $\zS_1 <\tau F^j T_{z_0} <\zS_2 $ and let $x=(x_0,i)$ be an admissible sink in $\zS_1^+$ such that \[ \zS_1^+ <\tau F^i T_{x_0} <\zS_2 .\]
We claim that $\dc(\sx\zS_1^+,\zS_2)<\dc(\zS_1,\zS_2)$.

We first prove that $G_x<\zS_2$ (see section \ref{sect 3.2} for the notation $G_x$). By definition, $G_x$ is constructed by taking closures under socle factors of injectives (lying on the slice) and predecessors. Taking predecessors (of  predecessors) of $\zS_2$ cannot create elements of $\zS_2$ or successors of $\zS_2$. Therefore, it suffices to show that, if $I$ is an injective predecessor of $\zS_2$ and $I\to M$, then $M<\zS_2$. Suppose that this is not the case, then $M\in \zS_2$ and, since $\zS_2$ is a local slice and $I$ is injective, then $I$ must belong to $\zS_2$, a contradiction.

Now the same argument as in case (1) above completes the proof of (a). Similarly, in case $\zS_2<\tau F^j T_{z_0} <\zS_1 $, we get (b).\qed
\end{pf}
\subsection{The main result}
We may now state and prove our main theorem.
\begin{thm}\label{thm 4.9}
Let $C$ be a tilted algebra having a tree $\zS$ as complete slice. The following conditions are equivalent:
\begin{itemize}
 \item[\textup{(a)}] $C'$ is a fibre quotient of $\Ctilde$.
\item[\textup{(b)}]  $C'$ is a fibre quotient of $\Crep$.
\item[\textup{(c)}] There exists a sequence of reflections and coreflections $\zs_1,\ldots,\zs_t$ such that $C'=\zs_1\cdots \zs_t C$ has $\zS'=\zs_1\cdots \zs_t\zS$ as complete slice and $C'=\Ctilde/\Ann\zS'$.
 \end{itemize}
\end{thm}
\begin{pf}
 Since the equivalence of (a) and (b) follows from Lemma \ref{lem 4.5}, and since Proposition \ref{prop 3.6} yields easily that (c) implies (a), it suffices to prove that (a) implies (c).
 
 Let $C'$ be a fibre quotient of $\Ctilde$. Then there exist two local slices $\zS$ and $\zS''$ in $\textup{mod}\,\Ctilde$ such that $C=\Ctilde/\Ann \zS$ and $C'=\Ctilde/\Ann \zS''$ (because of \cite{ABS2}). Lifting this information to $\Crep$, there exist two local slices $\stackrel{\scriptscriptstyle\vee}{\zS}$ and $\stackrel{\scriptscriptstyle\vee}{\zS''}$ lying in the same transjective component of $\zG(\textup{mod}\,\Crep)$ such that $G_\zl \stackrel{\scriptscriptstyle\vee}{\zS} =\zS$ and  $G_\zl \stackrel{\scriptscriptstyle\vee}{\zS''} =\zS''$. Applying Lemma \ref{lem 4.8} and an obvious induction, the finiteness of the distance function yields a sequence of reflections and coreflections $\zs_1,\ldots,\zs_t$ such that $\dc(\zs_1\cdots \zs_t \stackrel{\scriptscriptstyle\vee}{\zS} ,  \stackrel{\scriptscriptstyle\vee}{\zS''})=0$. This implies that $d(\zs_1\cdots \zs_t \zS,\zS'')=0.$ Let $\zS'=\zs_1\cdots \zs_t \zS$.
 By Proposition \ref{prop 3.6}, $C'=\zs_1\cdots \zs_t C$ is tilted and has $\zS'$ as a  complete slice. Let $C^*=\Ctilde/\Ann \zS'$, then $d(\zS',\zS'')=0$ implies $d(C^*,C')=0$. Because of Corollary \ref{cor 4.4}, we get indeed $C'=C^*$. This completes the proof.\qed
\end{pf}

\subsection{Example}\label{ex 4.10}
Let again $\Ctilde $ be the cluster-tilted algebra of Example \ref{ex 4.7}. We assume that $C$ is the tilted algebra given by the quiver 

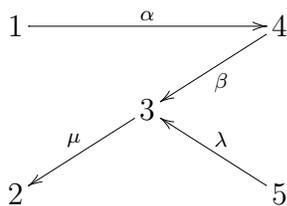
\[\xymatrix@R=20pt@C=40pt{
1\ar[rr]^\za&&4\ar[ld]^\zb\\
&3\ar[ld]_\mu\\
2&&5\ar[lu]_\zl
}
\]
bound by $\za\,\zb=0, \zl\,\mu=0$. A rightmost complete slice $\zS$ of $\textup{mod}\,C$ is given by 
 \[\zS=\left\{ \ {\begin{array}{c}4 \vspace{-10pt} \\ \vspace{-10pt} 3\\2\end{array}}\ ,\ 
{\begin{array}{c}4 \vspace{-10pt} \\3\end{array}}\ ,\ 
{\begin{array}{c}4\ 5 \vspace{-10pt} \\3\end{array}}\ ,\  
4\ ,\ 
{\begin{array}{c}1 \vspace{-10pt} \\4\end{array}}
\ \right\} \]
Reflecting successively at all admissible sinks yields successively the local slices
\[\begin{array}{rcl}
\zs_2\zS &=& \left\{ \ 
{\begin{array}{c}4\ 5 \vspace{-10pt} \\3\end{array}}\ ,\  
5\ , \ 4\ ,\ 
{\begin{array}{c}2 \vspace{-10pt} \\ 5\end{array}}\ ,\ 
{\begin{array}{c}1 \vspace{-10pt} \\4\end{array}}
\ \right\},
\\ 
\zs_3\zs_2\zS &=& \left\{ \ 
{\begin{array}{c}2 \vspace{-10pt} \\5\end{array}}\ ,\  
{\begin{array}{c}1 \vspace{-10pt} \\4\end{array}}
\ ,\ 
2\ , \ 1\ ,\ 
{\begin{array}{c}3 \vspace{-10pt} \\ 1\ 2 \end{array}}
\ \right\},
\\ 
\zs_4\zs_3\zs_2\zS &=& \left\{ \ 
{\begin{array}{c}2 \vspace{-10pt} \\5\end{array}}\ ,\  
2\ , \ 
{\begin{array}{c}3 \vspace{-10pt} \\ 1\ 2 \end{array}}\ ,\ 
{\begin{array}{c}3 \vspace{-10pt} \\2\end{array}}\ ,\ 
{\begin{array}{c}4 \vspace{-10pt} \\3 \vspace{-10pt} \\2\end{array}}
\ \right\},
\\ 
\zs_5\zs_3\zs_2\zS &=& \left\{ \ 
{\begin{array}{c}1 \vspace{-10pt} \\4\end{array}}\ ,\  
1\ , \ 
{\begin{array}{c}3 \vspace{-10pt} \\ 1\ 2 \end{array}}\ ,\ 
{\begin{array}{c}3 \vspace{-10pt} \\1\end{array}}\ ,\ 
{\begin{array}{c}5 \vspace{-10pt} \\3 \vspace{-10pt} \\1\end{array}}
\ \right\}.
\\ 
\zs_5\zs_4\zs_3\zs_2\zS &=& \left\{ \ 
{\begin{array}{c}3 \vspace{-10pt} \\ 1\ 2 \end{array}}\ ,\ 
{\begin{array}{c}3 \vspace{-10pt} \\1\end{array}}\ ,\ 
{\begin{array}{c}3 \vspace{-10pt} \\2\end{array}}\ ,\ 
{\begin{array}{c}5 \vspace{-10pt} \\3 \vspace{-10pt} \\1\end{array}}\ ,\ 
{\begin{array}{c}4 \vspace{-10pt} \\3 \vspace{-10pt} \\2\end{array}}
\ \right\}.
\\ 
\end{array}
\]
Then we have $\zS'=\zs_5\zs_4\zs_3\zs_2\zS=\zs_4\zs_5\zs_3\zs_2\zS$. The rightmost slice corresponding to $\zS'$ is 
\[ \zS'^+=\left\{ \ {\begin{array}{c}4 \vspace{-10pt} \\ \vspace{-10pt} 3\\2\end{array}}\ ,\ 
{\begin{array}{c}5 \vspace{-10pt} \\3 \vspace{-10pt} \\1\end{array}}\ ,\ 
{\begin{array}{c}4 \vspace{-10pt} \\3\end{array}}\ ,\ 
{\begin{array}{c}5 \vspace{-10pt} \\3\end{array}}\ , \ 
{\begin{array}{c}4\ 5 \vspace{-10pt} \\3\end{array}}  
\ \right\}, 
\]
therefore 
\[\zs_2 \zS'^+=\left\{ \
{\begin{array}{c}5 \vspace{-10pt} \\3 \vspace{-10pt} \\1\end{array}}\ ,\ 
{\begin{array}{c}5 \vspace{-10pt} \\3\end{array}}\ , \ 
{\begin{array}{c}4\ 5 \vspace{-10pt} \\3\end{array}}  \ ,\ 
5 \ ,\ 
{\begin{array}{c}2 \vspace{-10pt} \\5\end{array}}
\ \right\}, \ \ \ 
\]
while $\zs_1\zS'^+=\zS$. Therefore the fibre quotients of $\Ctilde $ are the algebras
\begin{enumerate}
\item $\zs_2 C$ given by the quiver       
\[
\xymatrix@R=20pt@C=40pt{
1\ar[rr]^\za&&4\ar[ld]^\zb\\
&3 \\
2\ar[rr]_\nu&&5\ar[lu]_\zl
}
\]
bound by  $\za\,\zb=0$   and $\nu\,\zl=0$.    
\item $\zs_3\zs_2C$ given by the quiver 
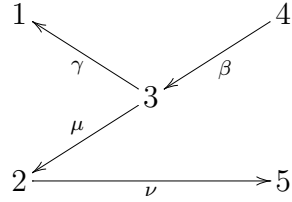
\[
\xymatrix@R=20pt@C=40pt{
1\ar[rr]^\za&&4 \\
&3\ar[lu]^\zg\ar[ld]_\mu\\
2\ar[rr]_\nu&&5 
}
\]
bound by    $\zg\,\za=0$  and $\mu\,\nu=0$. 
\item $\zs_4\zs_3\zs_2C$ given by the quiver 
\[
\xymatrix@R=20pt@C=40pt{
1 &&4\ar[ld]^\zb\\
&3\ar[lu]^\zg\ar[ld]_\mu\\
2\ar[rr]_\nu&&5 
}
\]
bound by    $\zb\,\zg=0$  and $\mu\,\nu=0$. 
\item  $\zs_5\zs_3\zs_2C$ given by the quiver 

\[
\xymatrix@R=20pt@C=40pt{
1\ar[rr]^\za&&4 \\
&3\ar[lu]^\zg\ar[ld]_\mu\\
2 &&5\ar[lu]_\zl
}
\]
bound by    $\zg\,\za=0$   and $\zl\,\mu=0$. 

\item  $\zs_5\zs_4\zs_3\zs_2C=\zs_4\zs_5\zs_3\zs_2C$ given by the quiver 
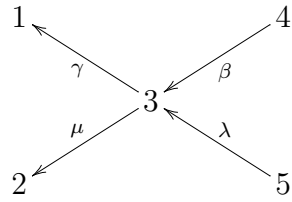
\[
\xymatrix@R=20pt@C=40pt{
1 &&4\ar[ld]^\zb\\
&3\ar[lu]^\zg\ar[ld]_\mu\\
2 &&5\ar[lu]_\zl
}
\]
bound by   $\zb\,\zg=0$   and $\zl\,\mu=0$. 

\item $\zs_2\zs_5\zs_4\zs_3\zs_2C$ given by the quiver 
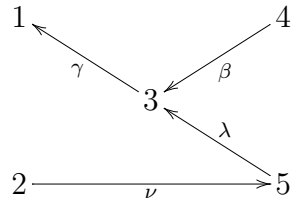
\[
\xymatrix@R=20pt@C=40pt{
1 &&4\ar[ld]^\zb\\
&3\ar[lu]^\zg \\
2\ar[rr]_\nu&&5\ar[lu]_\zl
}
\]
bound by $\zb\,\zg=0$ and $\nu\,\zl=0$. 
\end{enumerate}

Finally  $\zs_1\zs_5\zs_4\zs_3\zs_2C=C.$
  It is easily seen that we so obtain all fibre quotients of $\Ctilde$.

The reader can easily locate these reflections (fibre quotients) of $C$ in the quiver of $\Crep$:
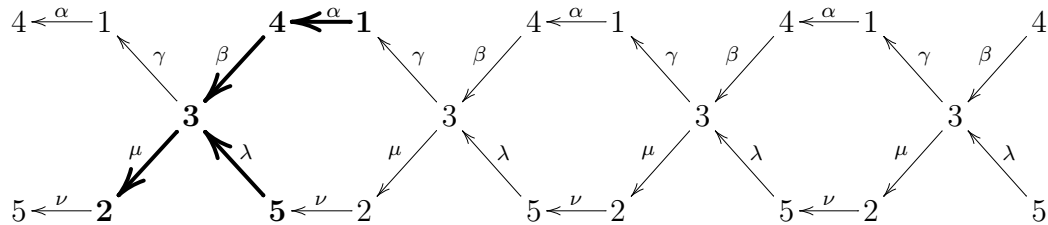
\[\xymatrix@C22pt{4&1\ar[l]_\za &&\mathbf{4}\ar@*{[|<1.5pt>]}[dl]_{{\zb}}&\mathbf{1}\ar@*{[|<1.5pt>]}[l]_(0.35){\za} &&4\ar[dl]_\zb&1\ar[l]_\za &&4\ar[dl]_\zb&1\ar[l]_\za &&4\ar[dl]_\zb\\
&&\mathbf{3}\ar[ul]_\zg\ar@*{[|<1.5pt>]}[dl]_\mu&&&3\ar[ul]_\zg\ar[dl]_\mu&&&3\ar[ul]_\zg\ar[dl]_\mu&&&3\ar[ul]_\zg\ar[dl]_\mu&&&\\
5&\mathbf{2}\ar[l]_\nu &&\mathbf{5}\ar@*{[|<1.5pt>]}[ul]_\zl&2\ar[l]_\nu&&5\ar [ul]_\zl&2\ar[l]_\nu&&5\ar[ul]_\zl&2\ar[l]_\nu&&5\ar[ul]_\zl}
\]
bound by the lifted relations $\za\,\zb=0$, $\zb\,\zg=0$, $\zg\,\za=0$ $\zl\,\mu=0$, $\mu\,\nu=0$ and $\nu\,\zl=0$. 
We have illustrated one copy of $C$ in bold face.


\section{Algorithm}\label{sect 5}
\subsection{} Let $C$ be a tilted algebra of tree type, and $\zG$ a connecting component of $\textup{mod}\,C$. We recall that a tilted algebra has a unique connecting component, except if it is concealed, in which case it has two. We denote by $\Sp$ and $\zS^-$, respectively, the rightmost and leftmost slice in $\zG$. We assume both $\Sp$ and $\zS^-$ exist. Let $\zG_1$ be the full subquiver of $\zG$ having as points
\[ \zG_1=\{M\in \ind C\mid \tau\zS^-\le M\le \tau^{-1}\Sp\}.\]

\begin{lem}
\label{lem 5.1}
 With the above notation,
 \begin{itemize} \item[\textup{(a)}] $\zG_1$ embeds as a full subquiver of $\zG(\textup{mod}\,\Crep)$.
   \item[\textup{(b)}] Let $M$ be a $\Crep$-module such that
   $\tau\zS^- \le M \le \tau^{-1}\Sp$
 then $M$ is a $C$-module lying in  $\zG_1$. 
  \end{itemize}
\end{lem}
\begin{pf}
(a) follows from Proposition \ref{prop 2.2}.

(b) Let $M$ be such a $\Crep$-module. It follows from the structure of $\zG(\textup{mod}\,\Crep)$ that $M$ lies in a transjective component and furthermore there exists $t\ge 0$ such that $\tau_{\Crep}^{-t} M\in \Sp$, that is, there exists a $C$-module $N\in \Sp$ such that $\tau_{\Crep}^{-t} M=N$. Applying Proposition \ref{prop 2.2}, we get $M= \tau_{\Crep}^{t}\tau_{\Crep}^{-t} M = \tau_{\Crep}^{t} N\cong  \tau_{C}^{t} N$, hence the statement.\qed
\end{pf}

\begin{rem}
 Note that if, for instance, $\zS^-$ does not exist, but $\Sp$ does, then the statement of the Lemma applies to the full subquiver of $\zG$ with points $\{M\in \ind C\mid M\le \tau^{-1}\Sp\}$.
\end{rem}
\subsection{}
Let now $x$ be an admissible sink in $C$ such that $G_x$ is contained in the rightmost slice $\Sp$ of  $\textup{mod}\,C$. Let $ I_y$ be a source in $G_x$ and define a $\Crep$-module $\Pbar_y $ by
\[\begin{array}{rcl}
\textup{top}\, \Pbar_y&=&S_y\\
\rad \Pbar_y &=& \tau^{-1}_C(I_y/S_y) = \displaystyle\bigoplus _{I_y\to M} (\tau^{-1}_CM).
\end{array}
\]
Note that, since $I_y$ is a source, then all indecomposable modules $M$ such that there exists an arrow $I_y\to M$ in $\zG(\textup{mod}\,C)$ lie in $G_x$ (see Section \ref{sect 3.2}). Also, as morphisms from $\textup{top}\, \Pbar_y$ to $\rad \Pbar_y$, we take, for every arrow $\za:y\to z$, the linear map $f_\za:\Pbar_y(y)\to \Pbar_y(z)$ defined by the right multiplication by the residual class of the arrow $\za$ in $\Crep=k\!\!\stackrel{\scriptscriptstyle\vee}{Q}\!/\!\stackrel{\scriptscriptstyle\vee}{I}$.

Recursively, for every $I_z$ in $G_x$ with the property that for each predecessor $I_w$ of $I_z$ in $G_x$, we have already introduced a corresponding projective module $\Pbar_w$, we  define $\Pbar_z$ by
\[\begin{array}{rcl}
\textup{top}\, \Pbar_z&=&S_z\\
\rad \Pbar_z &=& \tau^{-1}_C(I_z/S_z)\bigoplus\left( \displaystyle\bigoplus _{I_w\to I_z} \Pbar_w\right),
\end{array}
\]
where the second direct sum is taken over all arrows $I_w\to I_z$ in $G_x$.

Again, for  morphisms from $\textup{top}\, \Pbar_z$ to $\rad \Pbar_z$, we take, for every arrow $\za:z\to v$, the linear map $f_\za:\Pbar_z(z)\to \Pbar_z(v)$ defined by the right multiplication by the residual class of the arrow $\za$ in $\Crep=k\!\!\stackrel{\scriptscriptstyle\vee}{Q}\!/\!\stackrel{\scriptscriptstyle\vee}{I}$. The module $\Pbar_z$ is thus located at the position $\tau^{-2}I_z $ in $\zG(\textup{mod}\,\Crep)$.

\begin{lem}
\label{lem 5.2}
For each injective module $I_y$ in $G_x$, the $\Crep$-module $\Pbar_y$ thus constructed is isomorphic to the indecomposable projective $\Crep$-module $\Prep_y$ at $y$.
\end{lem}
\begin{pf}
 Clearly, it suffices to show that $\rad\Prep_y=\rad \Pbar_y$.  We have that $\rad\Prep_y$ is the direct sum of all $N\in \ind \Crep$ such that there exists an arrow $N\to \Prep_y$ in $\zG(\textup{mod}\,\Crep)$. There are two possibilities for such a radical summand $N$:

Either $N$ is not projective, and then there exists an arrow $I_y\to M$ with $M\cong \tau_{\Crep} N$ because $\Prep_y$ is also situated at the position $\tau^{-2} I_y$ in $\zG(\textup{mod}\,\Crep)$ (see Lemma \ref{lem 2.4}\,(a)),
or $N=\Prep_w$ is projective, and then there exists an arrow $\Prep_w\to\Prep_z$ in $\zG(\textup{mod}\,\Crep)$. 

Thus 
\[\rad \Prep_y =\left(\displaystyle\bigoplus_{I_y\to M} \tau^{-1}_{\Crep} M\right) \bigoplus \left(\bigoplus_{\Prep_w\to \Prep_z} \Prep_w\right),\]
where the two direct sums are taken over arrows in  $\zG(\textup{mod}\,\Crep)$.

Now, if $I_y=\Irep_y$ is a source in $G_x$, then there is no arrow $I_z\to I_y$ in $\zG(\textup{mod}\,C)$ and, because of Lemma \ref{lem 5.1}, there is no arrow $\Irep_z\to\Irep_y$ in $\zG(\textup{mod}\,\Crep)$. By Lemma \ref{lem 2.4}\,(b), there is no arrow $\Prep_z\to  \Prep_x$ in  $\zG(\textup{mod}\,\Crep)$. Therefore, using Proposition \ref{prop 2.2},
\[ \rad_{\Crep}  \Prep_y = \bigoplus_{\Irep_y \to M} \tau^{-1}_{\Crep} M = \bigoplus_{I_y\to M} \tau^{-1}_C M = \rad_{\Crep}\Pbar_y,\]
where the first direct sum is taken over arrows in $\zG(\textup{mod}\,\Crep) $ and the second over arrows in $\zG(\textup{mod}\,C)$.    

Now assume that $I_z$ is not a  source in $G_z$, By induction, we may suppose that $\Prep_w=\Pbar_w$ for all $w$ such that $I_w $ precedes $I_z$ in $G_x$. Thus 
\[\bigoplus_{\Prep_w\to\Prep_z} \Prep_w \cong \bigoplus_{\Irep_w\to\Irep_z} \Prep_w  \cong \bigoplus_{I_w\to I_z} \Prep_w \cong \bigoplus_{I_w\to I_z} \Pbar_w,
\]
where the last equality holds by induction. Since we have, as before,
\[\bigoplus_{\Irep_z\to M} \tau^{-1}_{\Crep} M = \bigoplus_{I_z\to M} \tau^{-1}_{C} M,\]
the proof is complete. \qed
 \end{pf}
\subsection{}
\begin{cor}
\label{cor 5.3} With the above notation, we have 
\[\sx\Sp = \{\zS\setminus G_x\}\cup\{\Pbar_y\mid I_y\in G_x \textup{ injective}\}\cup \{\tau^{-1}_C M\mid M\in G_x \textup{ not injective}\}.\]
\end{cor}
\begin{pf}
This follows directly from Lemma \ref{lem 5.2} and the construction in Section \ref{sect 3.3}.
\end{pf}
\begin{rem}
Clearly, the dual construction, starting from an admissible source $y$ in $C$ and constructing the local slice $\zs_y\zS^-$ in $\zG(\textup{mod}\,\Crep)$ holds as well. We leave its statement to the reader.
\end{rem}
\subsection{}
 We now describe an algorithm allowing to construct the transjective component of a cluster-tilted algebra $\Ctilde$ knowing only a complete slice of a tilted algebra $C$. Since the pushdown functor $G_\zl:\textup{mod}\,\Crep \to \textup{mod}\,\Ctilde$ is dense and thus induces an isomorphism of quivers $\zG(\textup{mod}\,\Ctilde)\cong \zG(\textup{mod}\,\Crep)/\mathbb{Z}$ (see \cite{ABS3}), it suffices to construct a transjective component of $\Crep$.
 
 Let $\zS$  be a complete slice in $\textup{mod}\,C$, then $\zS$ embeds as a local slice in a  transjective component $\zG$ of the cluster repetitive algebra $\Crep$. For clarity, we treat separately the representation-finite and the representation-infinite case.
 
\begin{itemize} \item[\textup{(a)}] \underline{Assume $\Ctilde$ is representation-finite}, that is, $\Crep$ is locally representation-finite. In this case, $\zS$ is a Dynkin quiver. We carry out the following steps.
\begin{enumerate}
\item If there exists a source $M$ of $\zS$ which is not injective, then we replace $\zS$ by 
\[\zS' =\{\zS\setminus\{M\}\}\cup \{\tau^{-1}M\}
\]
(here, the Auslander-Reiten translation $\tau $ is computed with respect to the support of $\zS$ which, at the start, is equal to $C$). 
If not go to 2. Repeat until  every source is injective.
\item If all sources of $\zS$ are injective then there exists a source $I_x$ in $\zS$ such that $G_x$ exists (because of Lemma \ref{lem 3.2}). Then we replace $\zS$ by
\[\zS' = \sx\zS.\]
Go to 1.
\end{enumerate}

 Since $\Crep$ is locally representation-finite, we   eventually construct a slice $\zS$ such that for every module $M$ in $\zS$, the module $\varphi^{-1} M$  has already been constructed before, where $\varphi$ is the automorphism of $\Crep$ inducing the covering $\Crep\to \Ctilde$ (see Section \ref{sect 2.3}). At this point the algorithm stops.  After identification under $\varphi$, we have thus obtained the Auslander-Reiten quiver of the cluster-tilted algebra $\Ctilde$. 

\item[\textup{(b)}] \underline{Assume $\Ctilde$ is representation-infinite}, that is, $\Crep$ is locally representation-infi\-nite.  We carry out the following steps.
\begin{enumerate}
\item If there exists a source $M$ of $\zS$ which is not injective, then we replace $\zS$ by 
\[\zS' =\{\zS\setminus\{M\}\}\cup \{\tau^{-1}M\}
\]
(where, again,  $\tau^{-1} $ is computed with respect to the support of $\zS$). 
Repeat. If this procedure produces a $\zS$ in which every source is injective, then go to 2.
If not, then this procedure produces the right stable part of $\zG$. Then go to 3.
\item If all sources of $\zS$ are injective then there exists a source $I_x$ in $\zS$ such that $G_x$ exists. Then we replace $\zS$ by
\[\zS' = \zs_x\zS.\]
Go to 1.
Since there are finitely many injectives in $\zG$ then, at some point, we get to 3.
\item Return to the initial slice $\zS$.
\item If there exists a sink $N$ of $\zS$ which is not projective, then we replace $\zS$ by 
\[\zS' =\{\zS\setminus\{N\}\}\cup \{\tau N\}
\]
(where, again,  $\tau $ is computed with respect to the support of $\zS$). 
Repeat. 
If this procedure produces a $\zS$ in which every sink is projective, then go to 5.
If not, then this procedure produces the left stable part of $\zG$. Then the algorithm stops.
\item If all sinks of $\zS$ are projective then there exists a sink $P_y$ in $\zS$ such that $G_y$ exists. Then we replace $\zS$ by
\[\zS' = \zs_y\zS.\]
Go to 4.
Since there are finitely many projectives in $\zG$ then, at some point, the algorithm stops.
\end{enumerate}
 \end{itemize} 
\begin{thm}\label{thm 5.4}
Let $C$ be a tilted algebra of tree type. Then the transjective component of $\zG(\textup{mod}\,\Ctilde)$ is constructed by the preceding algorithm. Moreover, if $C$ is of Dynkin type, then the algorithm yields $\zG(\textup{mod}\,\Ctilde)$.
\end{thm}
\begin{pf}
This follows from Corollary \ref{cor 5.3} and the density of the pushdown functor $G_\zl :\textup{mod}\,\Crep \to \textup{mod}\,\Ctilde$. \qed\end{pf}

\subsection{A representation-finite example }\label{ex 5.4}
Let $C$ be the tilted algebra of type $\mathbb{D}_4$ given by the quiver

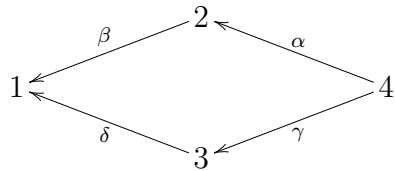
\[\xymatrix@C=60pt@R=15pt{
&2\ar[dl]_\zb\\
1 &&4\ar[ul]_\za\ar[dl]^\zg\\
&3 \ar[ul]^\zd
}
\]
bound by $\za\zb=\zg\zd$. We construct its Auslander-Reiten quiver until we reach its rightmost slice 
\[\Sp= \left\{\
{\begin{array}{c} 4 \vspace{-10 pt} \\2\ 3 \vspace{-10 pt} \\1 \end{array}}\ ,\ 
{\begin{array}{c} 4 \vspace{-10 pt} \\2\ 3 \end{array}} \ ,\ 
{\begin{array}{c} 4 \vspace{-10 pt} \\2 \end{array}}\ ,\ 
{\begin{array}{c} 4 \vspace{-10 pt} \\3 \end{array}} \ \right\}. 
\]
Since $\Sp$ has a unique source ${\begin{array}{c} 4 \vspace{-10 pt} \\2\ 3 \vspace{-10 pt} \\1 \end{array}}$, the corresponding sink $1$ is admissible and so we get 
\[\zs_1^+\Sp= \left\{\
{\begin{array}{c} 4 \vspace{-10 pt} \\2 \end{array}}\ ,\ 
{\begin{array}{c} 4 \vspace{-10 pt} \\3 \end{array}}\ ,\ 
4, \ , \ 
{\begin{array}{c} 1 \vspace{-10 pt} \\4 \end{array}} 
\ \right\}. 
\]
In  the next step we must move the points ${\begin{array}{c} 4 \vspace{-10 pt} \\2 \end{array}}$ and ${\begin{array}{c} 4 \vspace{-10 pt} \\3 \end{array}}$ simultaneously (because $G_2=G_3$), hence we get
\[\zs_2^+\zs_1^+\Sp=\zs_3^+\zs_1^+\Sp= \left\{\
{\begin{array}{c} 1 \vspace{-10 pt} \\4 \end{array}} \ ,\ 
1 \ , \ 
{\begin{array}{c} 2 \vspace{-10 pt} \\1 \end{array}}\ ,\ 
{\begin{array}{c} 3 \vspace{-10 pt} \\1 \end{array}}
\ \right\}. 
\]
A further reflection yields \[\zs_4^+\zs_2^+\zs_1^+\Sp= \left\{\
{\begin{array}{c} 2 \vspace{-10 pt} \\1 \end{array}}\ ,\ 
{\begin{array}{c} 3 \vspace{-10 pt} \\1 \end{array}}\ ,\ 
{\begin{array}{c} 2\ 3 \vspace{-10 pt} \\1 \end{array}}\ ,\ 
{\begin{array}{c} 4 \vspace{-10 pt} \\2\ 3 \vspace{-10 pt} \\1 \end{array}}
\ \right\},
\]
which is the leftmost slice $\zS^-$ in $\zG(\textup{mod}\,C)$.
The Auslander-Reiten quiver of $\tilde C$ is of the form shown in Figure \ref{fig arq2}.
\begin{figure}
\[
\xymatrix@C=15pt@R=0pt{
&{\begin{array}{c} 2 \vspace{-10 pt} \\1 \end{array}}\ar[dr]
&&{\begin{array}{c} 3 \end{array}}\ar[dr]
&&{\begin{array}{c} 4 \vspace{-10 pt} \\2 \end{array}}\ar[dr]
&&
&&{\begin{array}{c} 2 \vspace{-10 pt} \\1 \end{array}}\\
{\begin{array}{c} \underline{1} \end{array}}\ar[dr]\ar[ur]
&
&{\begin{array}{c} 2\  3 \vspace{-10 pt} \\1 \end{array}}\ar[dr]\ar[ur]\ar[r]
&{\begin{array}{c} 4 \vspace{-10 pt} \\2\ 3 \vspace{-10 pt} \\1 \end{array}}\ar[r]
&{\begin{array}{c} 4 \vspace{-10 pt} \\2\ 3 \end{array}}\ar[dr]\ar[ur]
&
&{\begin{array}{c} 4 \end{array}}\ar[r]
&{\begin{array}{c} 1 \vspace{-10 pt} \\4 \end{array}}\ar[r]
&{\begin{array}{c} \underline{1} \end{array}}\ar[ru]\ar[rd]
&\\
&{\begin{array}{c} 3 \vspace{-10 pt} \\1 \end{array}}\ar[ru]&&
{\begin{array}{c} 2\end{array}}\ar[ru]&&
{\begin{array}{c} 4 \vspace{-10 pt} \\3 \end{array}}\ar[ru]&&
&&
{\begin{array}{c} 3 \vspace{-10 pt} \\1 \end{array}}
}
\]
\caption{Auslander-Reiten quiver of Example \ref{ex 5.4}}\label{fig arq2}
\end{figure}

\subsection{A representation-infinite example}
Let $C$ be the tilted algebra of type $\tilde{\mathbb{D}}_4$ given by the quiver
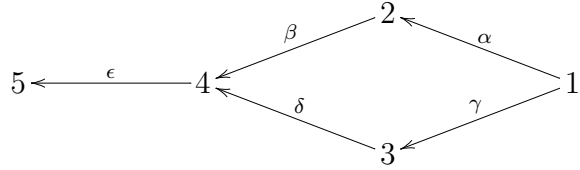
\[\xymatrix@C=60pt@R=15pt{&&2\ar[dl]_\zb\\
5&4\ar[l]_\ze &&1\ar[ul]_\za \ar[dl]_\zg\\
&&3\ar[ul]_\zd } \]
bound by $\za\,\zb =\zg\zd$ and $\za\,\zb\,\ze=0$.
Here, $\Ctilde $ is representation-infinite. We show part of its transjective component.
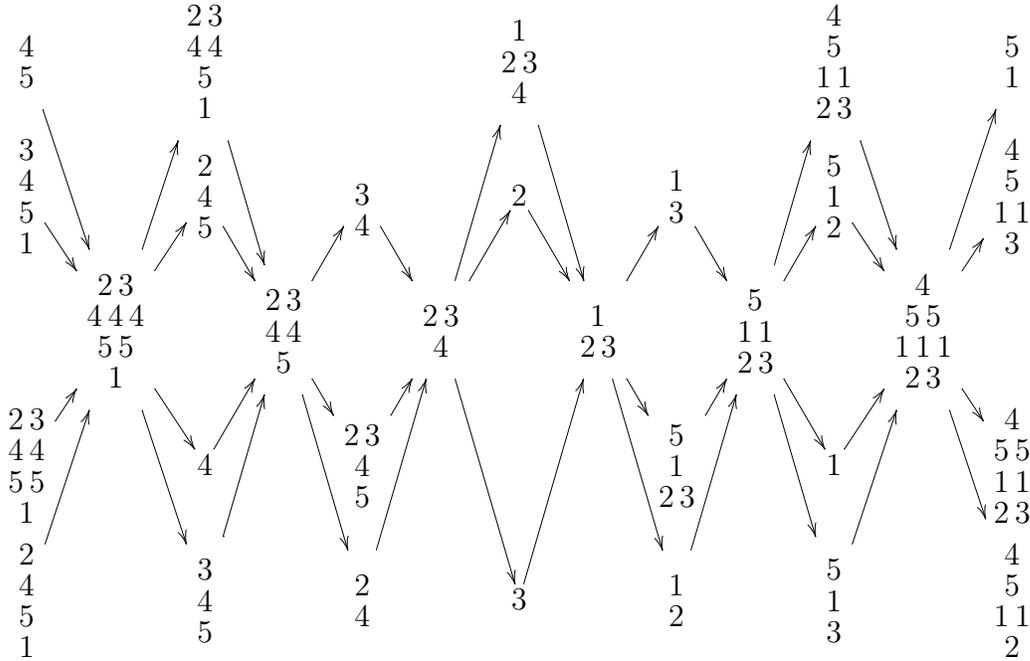
\[\xymatrix@C=8pt@R=-10pt{
{\begin{array}{c} 4 \vspace{-10 pt} \\5 \end{array}}\ar[ddr]
&&{\begin{array}{c} 2\,3\vspace{-10 pt}\\ 4\,4\vspace{-10pt}\\5 \vspace{-10 pt} \\1 \end{array}}\ar[ddr]
&&
&&{\begin{array}{c} 1\vspace{-10 pt}\\ 2\,3\vspace{-10pt}\\4 \end{array} }\ar[ddr]
&&
&&{\begin{array}{c} 4\vspace{-10 pt}\\ 5\vspace{-10 pt}\\ 1\,1\vspace{-10pt}\\2\,3 \end{array} }\ar[ddr]
&&{\begin{array}{c}  5\vspace{-10 pt}\\ 1 \end{array} }
\\
{\begin{array}{c} 3\vspace{-10 pt}\\ 4\vspace{-10 pt}\\ 5\vspace{-10pt}\\1 \end{array} }\ar[dr]
&&{\begin{array}{c} 2\vspace{-10 pt}\\ 4\vspace{-10 pt}\\ 5 \end{array} }\ar[dr]
&&{\begin{array}{c} 3\vspace{-10 pt}\\ 4\vspace{-10 pt} \end{array} }\ar[dr]
&& 2\ar[dr]
&&{\begin{array}{c} 1\vspace{-10 pt}\\ 3\end{array} }\ar[dr]
&&{\begin{array}{c} 5\vspace{-10 pt}\\ 1\vspace{-10 pt}\\ 2 \end{array} }\ar[dr]
&&{\begin{array}{c} 4\vspace{-10 pt}\\ 5\vspace{-10 pt}\\ 1\,1\vspace{-10pt}\\3 \end{array} }
\\
&{\begin{array}{c} 2\,3\vspace{-10 pt}\\ 4\,4\,4\vspace{-10 pt}\\ 5\,5\vspace{-10pt}\\ 1 \end{array} }\ar[uur]\ar[ur]\ar[dr]\ar[ddr]
&&{\begin{array}{c} 2\,3\vspace{-10 pt}\\ 4\,4\vspace{-10 pt}\\ 5 \end{array} }\ar[ur]\ar[dr]\ar[ddr]
&&{\begin{array}{c} 2\,3\vspace{-10 pt}\\ 4\end{array} }\ar[uur]\ar[ur]\ar[ddr]
&&{\begin{array}{c} 1\vspace{-10 pt}\\ 2\,3 \end{array} }\ar[ur]\ar[dr]\ar[ddr]
&&{\begin{array}{c} 5\vspace{-10 pt}\\ 1\,1\vspace{-10 pt}\\ 2\,3\end{array} }\ar[uur]\ar[ur]\ar[dr]\ar[ddr]
&&{\begin{array}{c} 4\vspace{-10 pt}\\ 5\, 5\vspace{-10 pt}\\ 1\,1\,1\vspace{-10pt}\\ 2\,3 \end{array} }\ar[uur]\ar[ur]\ar[dr]\ar[ddr]
\\
{\begin{array}{c} 2\,3\vspace{-10 pt}\\ 4\,4\vspace{-10 pt}\\ 5\,5\vspace{-10pt}\\ 1 \end{array} }\ar[ur]
&&4\ar[ur]
&&{\begin{array}{c} 2\,3\vspace{-10 pt}\\ 4\vspace{-10 pt}\\ 5\end{array} }\ar[ur]
&&
&&{\begin{array}{c} 5\vspace{-10 pt}\\ 1\vspace{-10 pt}\\ 2\,3 \end{array} }\ar[ur]
&&1\ar[ur]
&&{\begin{array}{c} 4\vspace{-10 pt}\\ 5\,5\vspace{-10 pt}\\ 1\,1\vspace{-10pt}\\ 2\,3\end{array} }
\\
{\begin{array}{c} 2\vspace{-10 pt}\\ 4\vspace{-10 pt}\\ 5\vspace{-10pt}\\ 1 \end{array} }\ar[uur]
&&{\begin{array}{c} 3\vspace{-10 pt}\\ 4\vspace{-10 pt}\\ 5 \end{array} }\ar[uur]
&&{\begin{array}{c} 2\vspace{-10 pt}\\4 \end{array} }\ar[uur]
&&3 \ar[uur]
&&{\begin{array}{c} 1\vspace{-10 pt}\\ 2 \end{array} }\ar[uur]
&&{\begin{array}{c} 5\vspace{-10 pt}\\ 1\vspace{-10 pt}\\3 \end{array} }\ar[uur]
&&{\begin{array}{c} 4\vspace{-10 pt}\\ 5\vspace{-10 pt}\\ 1\,1\vspace{-10pt}\\ 2 \end{array} }
}
\]
The rest of the transjective component  is constructed by the ``knitting'' procedure, constructing successively the Auslander-Reiten translates of the modules thus obtained. The remaining projectives lie in the tubes. The cluster repetitive algebra $\Crep$ is given by the quiver
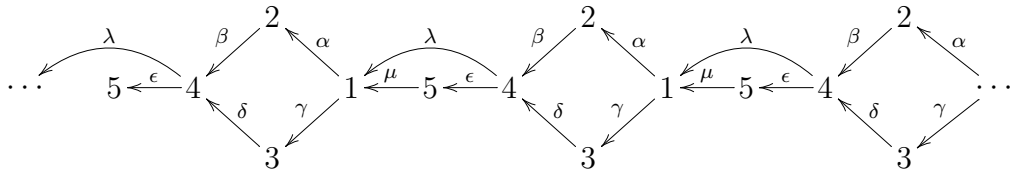
\[\xymatrix@R15pt@C20pt{&&&2\ar[dl]_\zb&&&&2\ar[dl]_\zb&&&&2\ar[dl]_\zb&&&&
\\
\cdots&5&4\ar[l]_\ze\ar@/_15pt/[ll]_\zl && 1\ar[dl]_\zg\ar[ul]_\za & 5\ar[l]_\mu &4\ar[l]_\ze\ar@/_15pt/[ll]_\zl && 1\ar[dl]_\zg\ar[ul]_\za & 5\ar[l]_\mu&4\ar[l]_\ze\ar@/_15pt/[ll]_\zl && \cdots\ar[dl]_\zg\ar[ul]_\za\\
&&&3\ar[ul]_\zd &&&&3\ar[ul]_\zd&&&&3\ar[ul]_\zd
}\]
bound by $\za\,\zb=\zg\,\zd,\  \za\,\zb\,\ze=0,\ \zb\,\zl=\zb\,\ze\,\mu,\ \zl\,\za=\ze\,\mu\,\za,\ \zd\,\zl=0$ and $\zl\,\zg=0$.

\section{Tubes}\label{sect 6}
The same algorithm seems to work for the tubes of the cluster-tilted algebras of Euclidean type. We have no proof of this fact but we give partial results and an example here.

Let $A$ be a hereditary algebra of Euclidean type and $T$ be a tilting $A$-module without preinjective summands. Assume that $T_i$ is a summand of $T$ that lies in a tube and such that $i$ is a source of $C=\End\!_AT$. Denote by $r$ the quasi length of $T_i$ and let $M$ be the quasi simple module that lies on the same ray as $T_i$ on the mouth of the tube.

\begin{lem}\label{tube1}
The immediate predecessor of $T_i$ on the semi-ray ending at $T_i$ is a summand of $T$.
\end{lem}
\begin{pf}
If $r=1$, then $M=T_i$ and the result holds since there is no such predecessor. If $r>1$, it follows from the assumption that $i$ is a source in $C$.\qed
\end{pf}

We denote this predecessor by $T_j$. Thus there is a sectional path $M\to\cdot\to \cdots\to T_j\to T_i$ of length $r-1$, and $M$ lies on the mouth of the tube.

\begin{lem}\label{tube2} In the above situation, we have
\[\Hom_A(T,\tau^2T_i)\cong\Hom_A(T,\tau^2M).\]
\end{lem}
\begin{pf}
Applying the functor $\Hom_A(T,-)$ to the short exact sequence \[0\to \tau^2 M\to \tau^2 T_i\to \tau T_j\to 0, \] 
the result follows from $\Hom_A(T,\tau T_j)=D\Ext_A(T_j,T)=0$.
\qed\end{pf}

\begin{lem}\label{tube3} In the above situation, let $\tilde I_i$ denote the indecomposable injective and $\tilde S_i$ the indecomposable simple module of the cluster-tilted algebra $C\ltimes \DCC$ corresponding to the point $i$. Then
\[\tilde I_i/\tilde S_i = \Hom_A(T,\tau^2 T_i).\]
\end{lem}
\begin{pf}
A straightforward computation shows that
\[\begin{array}{rcl}
\tilde I_i& =&\Hom_{\mathcal{C}}(T,\tau^2 T_i) \\
&=& \Hom_A(T,\tau^2 T_i)\oplus \Hom_{\mathcal{D}^b(\textup{mod}\,A)}(\tau T[-1],\tau^2 T_i)\\
&=& \Hom_A(T,\tau^2 T_i)\oplus D\Hom_A(\tau^2 T_i,\tau^2 T).
\end{array}
\]
The simple socle of $\tilde I_i$ corresponds in this description to the direct summand $D\Hom_A(\tau^2 T_i,\tau^2 T_i)$ of the second term. Thus 
\[\tilde I_i/\tilde S_i =  \Hom_A(T,\tau^2 T_i)\oplus D\Hom_A(\tau^2 T_i,\tau^2 \overline T),\]
where $\overline T\oplus T_i=T$. The statement   now follows, because  $\Hom_A(\tau^2 T_i,\tau^2 \overline T)=\Hom_A(T_i,\overline T) =0$, because $i$ is a source in $C$. 
\qed
\end{pf}

Now consider the image of the tube in the module category of the tilted algebra $C=\End\!_AT$. The $A$-modules $T_j$ and $T_i$ correspond to the indecomposable projective $C$-modules $P_j$ and $P_i$ respectively. Moreover $P_j$ is a direct summand of the radical of $P_i$. Since $P_i$ lies in a tube its radical $\rad P_i=P_j\oplus N$, for some indecomposable $C$-module $N$.
Since $i$ is a source, it follows from the construction of the tube in $\textup{mod}\,C$ from the tube in $\textup{mod}\,A$  that $\tau_C N=\Hom_A(T,\tau^2 M).$ 
\begin{lem}\label{tube4}  With the notation above,
\[\tilde I_i/\tilde S_i =\tau_C N.\]
\end{lem}
\begin{pf} 
 $\tau_C N=\Hom_A(T.\tau^2 M) = \Hom_A(T,\tau^2 T_i) = \tilde I_i/\tilde S_i$, where the second equality follows from Lemma \ref{tube2} and the last   from Lemma \ref{tube3}.
\qed
\end{pf}

This shows that at least in certain cases, a similar algorithm as for the transjective component can be used to construct the tubes of the cluster-tilted algebra. Starting from the tube of the tilted algebra, we use knitting to the left until we reach an indecomposable projective $C$-module $P_i$. We insert a new injective at the position $\tau^2 P_i$  by requiring that its socle quotient is equal to $\tau_C$ of the unique non-projective indecomposable summand of the radical of $P_i$ in mod $C$. Lemma \ref{tube4} shows that this module is actually the indecomposable injective module $\tilde I_i$ of the cluster-tilted algebra.

The arguments above will stop functioning if we come to another projective $P_\ell$ inside the same tube for which there is no sectional path from $P_\ell$ to $P_i$. The algorithm still seems to work in all the examples we have computed, but we do not know how to prove it.

\begin{exmp}
We conclude with an example of a tube.
Let $C$ be given by the quiver 
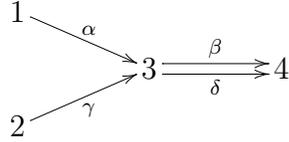
\[\xymatrix@R10pt@C40pt{1\ar[dr]^\za\\&3\ar@<2pt>[r]^\zb \ar@<-2pt>[r]_\zd&4 \\ 2\ar[ru]_\zg}\]
bound by the relations $\za\zb=0$ and $\zg\zd=0$.
One of the two exceptional tubes in $\textup{mod}\,C$ is given as
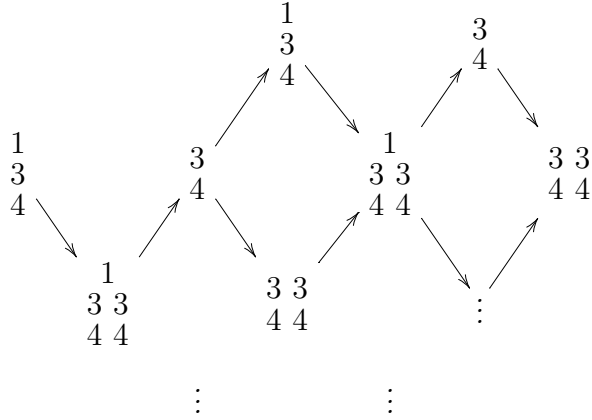
\[\xymatrix@R0pt@C15pt{ 
 &&&{\begin{array}{c} 1 \vspace{-10 pt} \\3\vspace{-10pt}\\4 \end{array}}\ar[rd] 
&&{\begin{array}{c} 3\vspace{-10pt}\\4 \end{array}} \ar[rd]
\\
 {\begin{array}{c} 1 \vspace{-10 pt} \\3\vspace{-10pt}\\4 \end{array}}\ar[rd] 
&& {\begin{array}{c} 3\vspace{-10pt}\\4 \end{array}} \ar[ru]\ar[rd]&&
 {\begin{array}{c} 1 \vspace{-10 pt} \\3\ 3\vspace{-10pt}\\4\ 4 \end{array}} \ar[rd]\ar[ru]&&
 {\begin{array}{c} 3\ 3\vspace{-10pt}\\4\ 4 \end{array}} 
\\
 &{\begin{array}{c} 1 \vspace{-10 pt} \\3\ 3\vspace{-10pt}\\4\ 4 \end{array}}\ar[ru] 
 &&{\begin{array}{c} 3\ 3\vspace{-10pt}\\4\ 4 \end{array}} \ar[ru]&&\vdots\ar[ru]\\
 &&\vdots&&\vdots
} \]
where modules with identical labels must be identified.
The module $P_1= {\begin{array}{c} 1 \vspace{-10 pt} \\3\vspace{-10pt}\\4 \end{array}}$ is projective and each module in the  tube lies in the $\tau$-orbit of $P_1$.

We use our algorithm to construct the tube of the corresponding cluster-tilted algebra $\tilde C=C\ltimes\DCC$ which is given by the quiver 
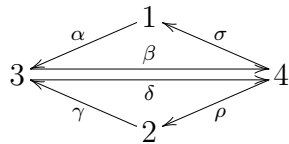
\[\xymatrix@R10pt@C40pt{&1\ar[dl]_\za\\3\ar@<2pt>[rr]^\zb \ar@<-2pt>[rr]_\zd&&4 \ar[lu]_\zs\ar[ld]^\rho\\& 2\ar[lu]^\zg}\]
bound by the relations $\za\zb=\zb\zs=\zs\za=\zg\zd=\zd\rho=\rho\zg=0$.
First we construct the new injective module
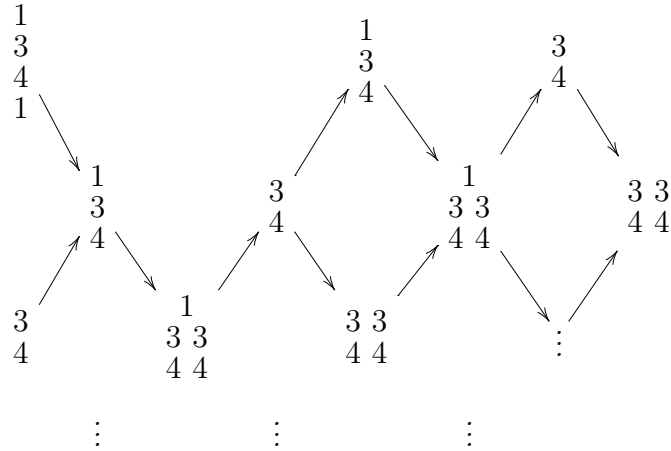
\[\xymatrix@R0pt@C15pt{ 
&{\begin{array}{c} 1 \vspace{-10 pt} \\3\vspace{-10pt}\\4 \vspace{-10pt} \\1\end{array}}\ar[rd] 
&&&&{\begin{array}{c} 1 \vspace{-10 pt} \\3\vspace{-10pt}\\4 \end{array}}\ar[rd] 
&&{\begin{array}{c} 3\vspace{-10pt}\\4 \end{array}} \ar[rd]
\\
&&{\begin{array}{c} 1 \vspace{-10 pt} \\3\vspace{-10pt}\\4 \end{array}}\ar[rd] 
&&{\begin{array}{c} 3\vspace{-10pt}\\4 \end{array}} \ar[ru]\ar[rd]&&
 {\begin{array}{c} 1 \vspace{-10 pt} \\3\ 3\vspace{-10pt}\\4\ 4 \end{array}} \ar[rd]\ar[ru]&&
 {\begin{array}{c} 3\ 3\vspace{-10pt}\\4\ 4 \end{array}} 
\\
& {\begin{array}{c} 3\vspace{-10pt}\\4 \end{array}} \ar[ru]
&&
 {\begin{array}{c} 1 \vspace{-10 pt} \\3\ 3\vspace{-10pt}\\4\ 4 \end{array}} \ar[ru]&&
 {\begin{array}{c} 3\ 3\vspace{-10pt}\\4\ 4 \end{array}} \ar[ru]&&\vdots\ar[ru]\\
 &&\vdots&&\vdots&&\vdots
} \]

and then we continue knitting to the left until the modules start repeating.

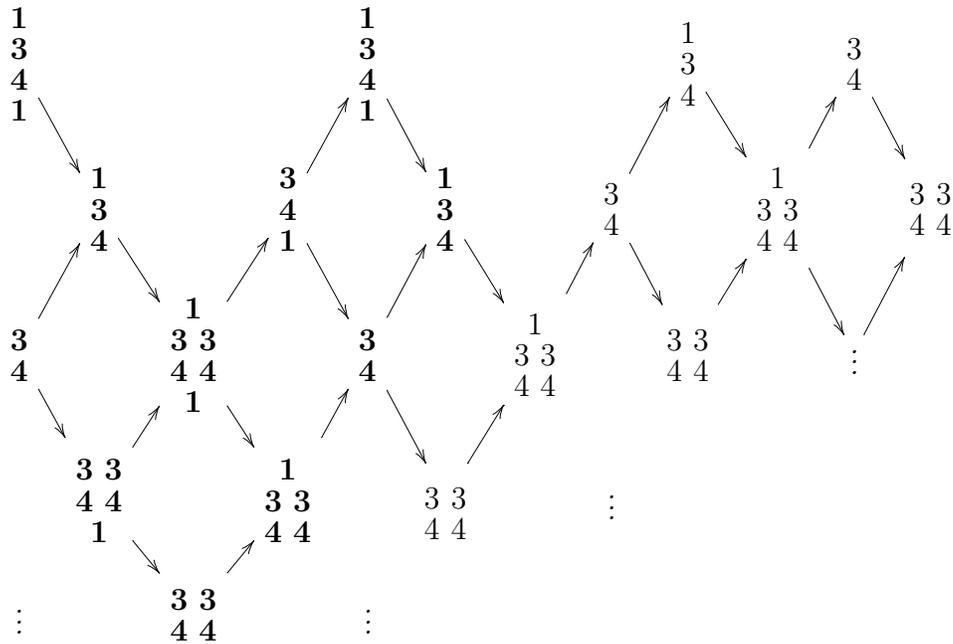
\[\xymatrix@R0pt@C10pt{ 
&{\begin{array}{c} {\bf 1} \vspace{-10 pt} \\{\bf 3}\vspace{-10pt}\\{\bf 4} \vspace{-10pt} \\{\bf 1}\end{array}}\ar[rd] 
&&&&{\begin{array}{c} {\bf 1} \vspace{-10 pt} \\{\bf 3}\vspace{-10pt}\\{\bf 4} \vspace{-10pt} \\{\bf 1}\end{array}}\ar[rd] 
&&&&{\begin{array}{c} 1 \vspace{-10 pt} \\3\vspace{-10pt}\\4 \end{array}}\ar[rd] 
&&{\begin{array}{c} 3\vspace{-10pt}\\4 \end{array}} \ar[rd]
\\
&& {\begin{array}{c} {\bf 1}\vspace{-10pt}\\{\bf 3}\vspace{-10pt}\\{\bf 4} \end{array}} \ar[rd]&&
{\begin{array}{c} {\bf 3} \vspace{-10 pt} \\{\bf 4}\vspace{-10pt}\\{\bf 1} \end{array}}\ar[rd]\ar[ru] 
&&{\begin{array}{c} {\bf 1} \vspace{-10 pt} \\{\bf 3}\vspace{-10pt}\\{\bf 4} \end{array}}\ar[rd] 
&&{\begin{array}{c} 3\vspace{-10pt}\\4 \end{array}} \ar[ru]\ar[rd]&&
 {\begin{array}{c} {1} \vspace{-10 pt} \\3\ 3\vspace{-10pt}\\4\ 4 \end{array}} \ar[rd]\ar[ru]&&
 {\begin{array}{c} 3\ 3\vspace{-10pt}\\4\ 4 \end{array}} 
\\
& {\begin{array}{c} {\bf 3}\vspace{-10pt}\\{\bf 4} \end{array}} \ar[ru]\ar[rd]
&& {\begin{array}{c} {\bf 1}\vspace{-10pt}\\{\bf 3}\ {\bf 3}\vspace{-10pt}\\{\bf 4}\ {\bf 4}\vspace{-10pt}\\{\bf 1} \end{array}} \ar[ru]\ar[rd]
&&{\begin{array}{c}  {\bf 3}\vspace{-10pt}\\{\bf 4} \end{array}}\ar[ru] \ar[rd]
&&
 {\begin{array}{c} 1 \vspace{-10 pt} \\3\ 3\vspace{-10pt}\\4\ 4 \end{array}}  \ar[ru]&&
 {\begin{array}{c} 3\ 3\vspace{-10pt}\\4\ 4 \end{array}} \ar[ru]&&\vdots\ar[ru]
 \\
&&
 {\begin{array}{c} {\bf 3}\ {\bf 3}\vspace{-10pt}\\{\bf 4}\ {\bf 4}\vspace{-10 pt} \\{\bf 1} \end{array}} \ar[ru]\ar[rd]&& {\begin{array}{c} {\bf 1}\vspace{-10pt}\\{\bf 3}\ {\bf 3}\vspace{-10pt}\\{\bf 4}\ {\bf 4} \end{array}} \ar[ru]
&& {\begin{array}{c} 3\ 3\vspace{-10pt}\\4\ 4\vspace{-10pt}\end{array}} \ar[ru]
&&\vdots
\\
&\vdots&& {\begin{array}{c}  {\bf 3}\ {\bf 3}\vspace{-10pt}\\{\bf 4}\ {\bf 4} \end{array}} \ar[ru]&&\vdots} \]

The tube in the cluster-tilted algebra consists of the modules in bold face. Modules (in bold face) with identical labels must be identified. Note that the  tube of the cluster-tilted algebra in this example is obtained by inserting a coray into the tilted tube.
\end{exmp}

\end{document}